\documentclass[11pt, a4paper]{article}
\usepackage{graphicx} 
\usepackage{hyperref}

\usepackage{eurosym}
\usepackage{amsfonts}
\usepackage{amsmath}
\usepackage[left=2cm, right=2cm, top=2cm, bottom=2cm]{geometry}
\usepackage{amssymb}
\usepackage{mathrsfs} 
\usepackage{physics}
\usepackage{siunitx}
\usepackage[dvipsnames]{xcolor}
\usepackage{tikz}
\usepackage{float}
\usepackage[affil-it]{authblk} 

\usepackage{amsthm}
\theoremstyle{plain}

\usepackage{caption}
\newcommand{\fr}{\text{Fr}}
\newcommand{\ro}{\text{Ro}}

\newcommand{\id}{\,\mathrm{d}}
\newcommand{\sdiff}{\,\circ\!\mathrm{d}}

\usepackage{titling}
\usepackage{graphicx}
\usepackage{duckuments}
\usepackage[format=plain, indention=0pt, width=0.85\textwidth, labelsep=colon,font=small,labelfont={bf}]{caption}
\usepackage{hyperref}
    \hypersetup{
        colorlinks=true,
        citecolor=Green,
        linkcolor=black,
        filecolor=magenta,      
        urlcolor=black,
        }
\usepackage{algorithm}
\usepackage{algorithmic}

\usepackage{subcaption}
\usepackage{booktabs, multirow}
\usepackage{colortbl}

\usepackage{scalerel}
\usepackage{pict2e}
\usepackage{tkz-euclide}
\usepackage{anyfontsize}
\usetikzlibrary{calc}
\usetikzlibrary{patterns, arrows.meta, decorations.pathreplacing}
\usetikzlibrary{shadows}
\usetikzlibrary{external}

\usepackage{pgfplots}
\pgfplotsset{compat=newest}
\usepgfplotslibrary{statistics}
\usepgfplotslibrary{fillbetween}

\definecolor{WaterBlue}{rgb}{0.11, 0.01 , 0.88}
\usepackage[dvipsnames]{xcolor}

\graphicspath{ {./figures/} }

\usepackage{lipsum}

\title{Bayesian inference for geophysical fluid dynamics using generative models}

\author{Alexander Lobbe\thanks{Department of Mathematics, Imperial College London} \qquad Dan Crisan\thanks{Department of Mathematics, Imperial College London}  \qquad Oana Lang\thanks{Department of Mathematics, Babes-Bolyai University} }

\date{}
\begin{document}

\maketitle

\abstract{Data assimilation plays a crucial role in numerical modeling, enabling the integration of real-world observations into mathematical models to enhance the accuracy and predictive capabilities of simulations. This approach is widely applied in fields such as meteorology, oceanography, and environmental science, where the dynamic nature of systems demands continuous updates to model states. However, the calibration of models in these high-dimensional, nonlinear systems poses significant challenges. In this paper, we explore a novel calibration methodology using diffusion generative models. We generate synthetic data that statistically aligns with a given set of observations (in this case the increments of the numerical approximation of a solution of a partial differential equation). This allows us to efficiently implement a model reduction and assimilate data from a reference system state modeled by a highly resolved numerical solution of the rotating shallow water equation of order $10^4$ degrees of freedom into a stochastic system having two orders of magnitude less degrees of freedom. To do so, the new samples are incorporated into a particle filtering methodology augmented with tempering and jittering for dynamic state estimation, a method particularly suited for handling complex and multimodal distributions. This work demonstrates how generative models can be used to improve the predictive accuracy for particle filters, providing a more computationally efficient solution for data assimilation and model calibration.
}

\section{Introduction}

Data assimilation is a technique in numerical modeling, in which real-world data is integrated into mathematical models to improve the accuracy and the predictive capabilities of simulations. This process has a wide range of applications in fields like meteorology, oceanography, or environmental science, where the dynamic nature of systems requires continual updates to model states based on incoming observations. However, the complexity of these systems often leads to challenges in accurately calibrating models, especially when dealing with sparse or noisy data.

Traditionally, model calibration involves adjusting model parameters to align with observational data, a process that can be computationally intensive and complex due to the high dimensionality and the nonlinear character of many systems. In response to these challenges, recent advancements in machine learning, particularly in the realm of generative models, seem to offer promising solutions. Diffusion generative models, known for their capacity to capture complex data distributions and generate realistic samples, present a novel approach for model calibration.

Our paper explores a novel methodology where a diffusion generative model is employed to calibrate a mathematical model prior to the commencement of the data assimilation process. By leveraging the diffusion model's ability to generate synthetic data that is statistically consistent with observed samples from high-dimensional phenomena, we achieve a more accurate and robust initial calibration of the mathematical model. This enhances the subsequent data assimilation process, leading to improved model performance and predictive accuracy.

\emph{Generative models} are a class of
machine learning models designed to generate new data samples from an unknown distribution which is typically available only through a dataset of samples. An important class of generative models are \emph{diffusion models} which gradually transform the training data into samples from a well-known distribution, such as a Gaussian distribution, using a mechanism analogous to diffusion. In this process, the forward and the
backward dynamics are learnt using a neural network. Once the learning is complete, samples from the unknown distribution are obtained by running the backward diffusion initiated from samples from the Gaussian distribution.

In the following, we will use the language of stochastic nonlinear filtering to explain the data assimilation methodology in general, and how it applies to the application discussed in this paper in particular. Therefore, we let $X$ and $Z$ be two processes defined on the probability space $(\Omega, \mathcal{F}, \mathbb{P})$. The process $X$ is usually called the \textit{signal process} or the \textit{truth} and $Z$ is the \textit{observation process}. In this paper, $X$ is the pathwise solution of a rotating shallow water system \eqref{rsw} approximated using a high-resolution numerical method.  The pair of processes $(X,Z)$ forms the basis of the nonlinear filtering problem which consists in finding the best approximation of the posterior distribution of the signal $X_t$ given the observations $Z_1, Z_2, \ldots, Z_t$ \footnote{For a mathematical introduction on the subject, see e.g. \cite{CrisanBain}. For an introduction from the data assimilation perspective of the filtering problem, see \cite{Leeuwenbook2015} and \cite{ReichCotter2015}.}. The posterior distribution of the signal at time $t$ is denoted by $\pi_t$. 

We let $d_X$ be the dimension of the state space and $d_Z$ be the dimension of the observation space.  In many real-world applications, particularly those related to weather prediction, $d_X$ 
is very large $d_X=O(10^9)$. Performing DA on such large dimensional models requires the use of super-computers and it is not surprising that some of the most advanced super-computers are used in meteorological offices around the world. 
Here we are advocating a different approach. Rather than working with the full signal $X_t$, we introduce an approximate model $X^c_t$ that is computed on a much coarser grid. In the example below, the signal is denoted by $X_t^f$ to emphasize the fact that it evolves on a finer grid than its proxy $X^c_t$ which is constructed on a coarser grid. Of course, $X_t^f$ and $X^c_t$ will behave quite differently as the small scale effects will be lost on the coarser scale. This is where generative modeling comes into play. To account for the effect of the small scale we add a stochastic term to the equation satisfied by $X^c_t$. The stochastic term will need to be calibrated to data obtained be recording $X_t^f$. The calibration needs to be done \emph{before} the data assimilation is applied. The two figures below illustrate this fact: 

\begin{figure}[h!]
    \centering
    \includegraphics[width=8cm, height=3.5cm]{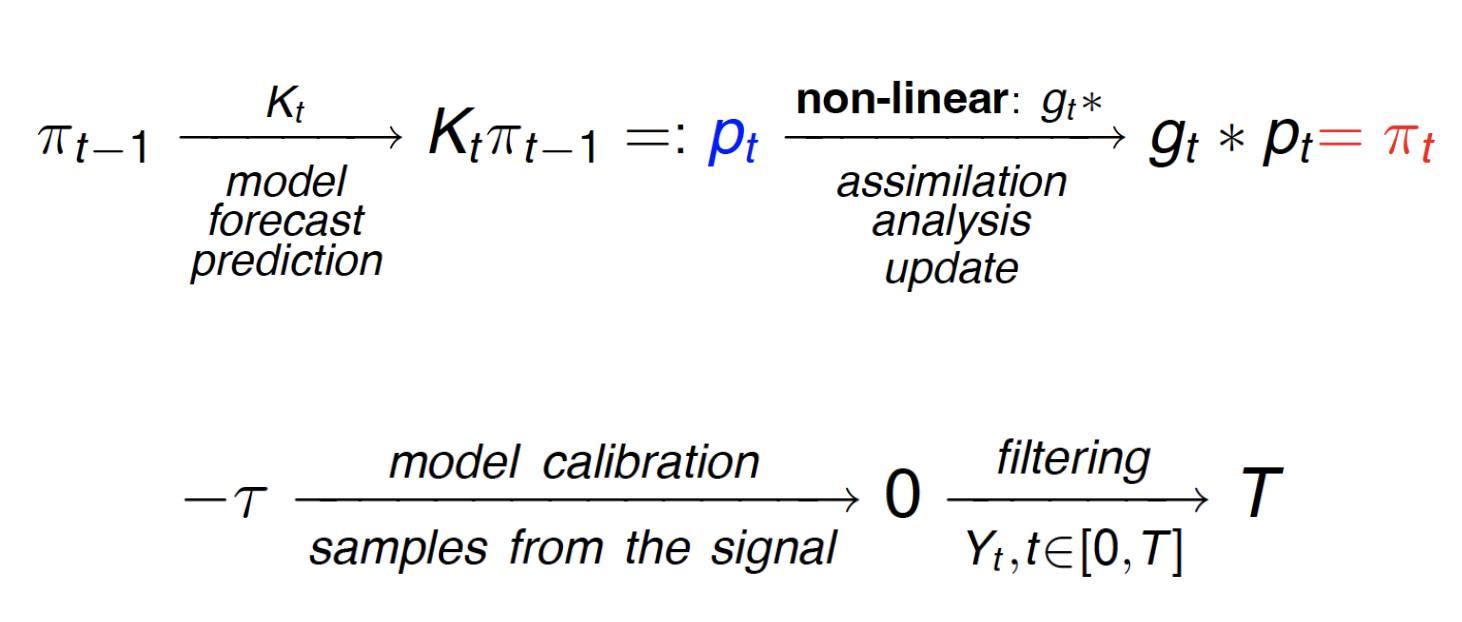}
    \caption{\textbf{Filtering and Calibration.} Data assimilation and model calibration. \emph{Top:} The predictive distribution is advanced using a forecast model and then updated via a nonlinear operation based on observed data. \emph{Bottom:} Before time $0$ the forecast model is calibrated. Filtering begins strictly after the calibration process is finshed.}
    \label{fig:filt-and-calib}
\end{figure}

We need to make an important remark here: The calibration methodology requires full and accurate knowledge of the signal, more precisely its projection on the coarse grid. For synthetic data this can be easily achieved by running the model equation on a sufficiently large time interval and on the fine grid and then projecting it on the coarser grid. For real data, one would need to use the so-called re-analysis data. Such data is freely available, for example, from 
ECMWF, e.g. ERA5 data (see e.g. \cite{era5}).

In \cite{Doucet} a generative model approach is used to approximate the posterior distribution directly, but in a procedure which must be performed offline. Repeating this iteratively would involve huge computational costs. In this work, we do not replace the two data assimilation steps (forecast and assimilation/update) using a generative model methodology, but rather use a diffusion model methodology to calibrate the signal \emph{before} starting the data assimilation process.

Let us explain briefly the methodology for the incorporation of the incomplete data. We will use particle filters. 
These are a class of filtering techniques which are adapted particularly to handle complex and multimodal distributions. Unlike traditional approaches, such as the Kalman filter, which assume linearity of the model and Gaussian noise, particle filters use a set of random samples/scenarios, or \emph{particles} to represent the posterior distribution of a system's state. Each particle carries a weight that reflects how well it fits the observed data. As new observations are assimilated, these particles are propagated through the model, and their weights are updated according to the \emph{likelihood} of the observed data given the particle's state: the particles that better match the observations are given higher weights, while those that are less accurate are down-weighted or discarded. Further details regarding this methodology can be found in \cite{LeeuwenCrisanLangPotthast} and in section \ref{sect:pf} below.

\begin{figure}[H]
    \centering
    \includegraphics[width=0.6\linewidth]{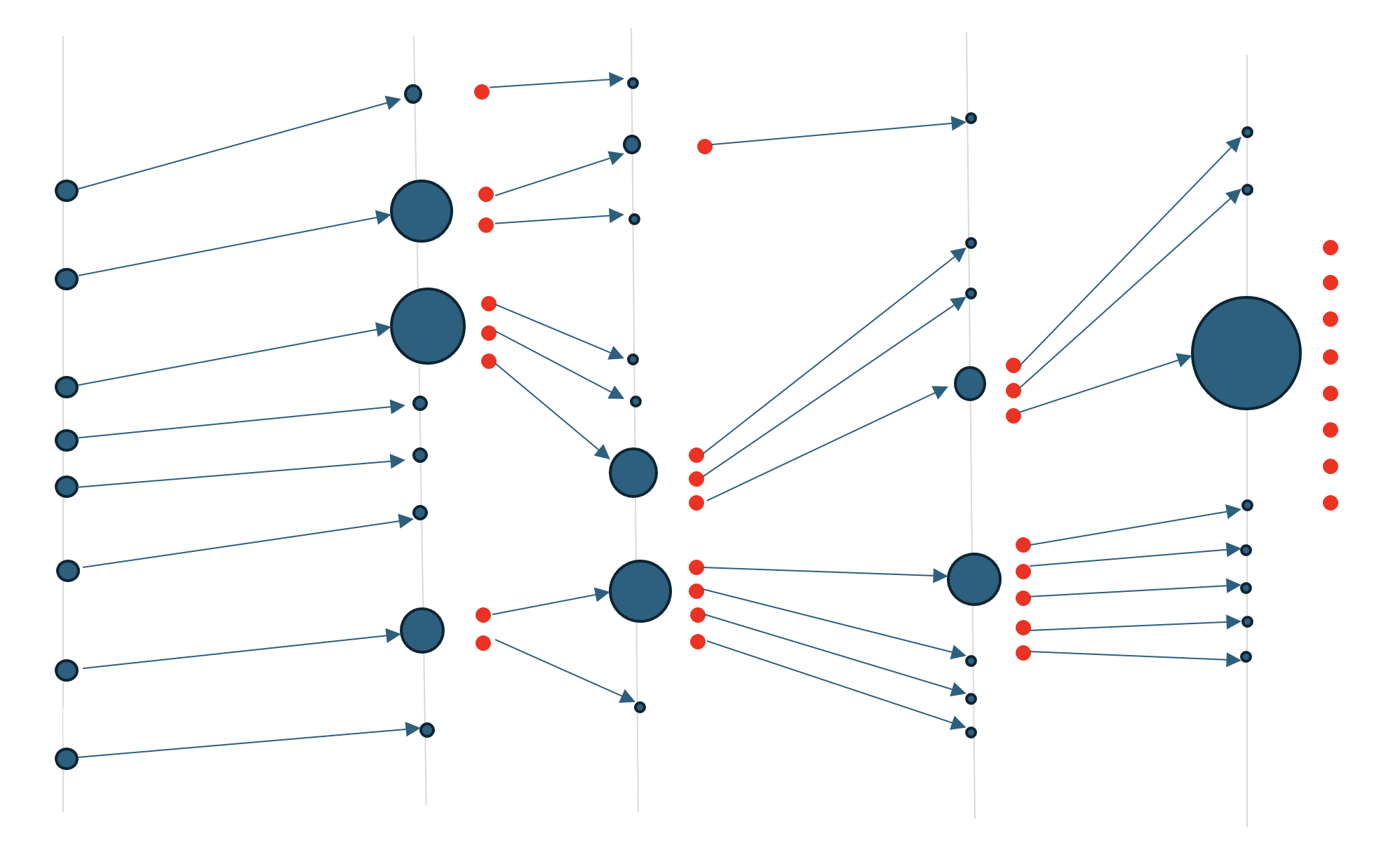}
   \caption{\textbf{Classical Particle Filter.} The resampling procedure ensures that particles with low weights are replaced with particles with higher weights. Following the resampling, an ensemble of equal-weighted particles is obtained. In high-dimensional spaces: one particle gains a weight close to one while all the others have weights close to zero and are discarded.}
 \label{classicalpf}
\end{figure}

\section{A Hybrid Rotating Shallow Water Model}

In this study, we develop our analysis based on a stochastic approximation of the following nondimensionalised rotating shallow water model:

\begin{equation}\label{rsw}
    \begin{aligned}
    & \id_t\vb{u} + \mathcal{R}(\vb{u},\eta )= 0\\
    & \id_t\eta + \mathcal{P}(\eta, \vb{u}) =0
    \end{aligned}
\end{equation}
where
\begin{itemize}
    \item $\vb{u}(x,t) = (u(x,t), v(x,t))$ is the horizontal fluid velocity vector field
    \item $\eta(x,t)$ is the height of the fluid column\footnote{We follow here the notation used in \cite{MPEbook}. In the previous work \cite{us}, the height of the fluid column was denoted by $h$.} 
\item The operator $\mathcal{R}: C^1(\Omega, \mathbb{R}^2) \times C^1(\Omega, \mathbb{R}) \to C^0(\Omega, \mathbb{R}^2)$, where $\Omega$ is a spatial domain, governs the dynamics of the velocity vector field:
$$\mathcal{R}(\vb{u}, \eta) := (\vb{u}\cdot\grad) \vb{u} + \frac{f}{\ro} \hat{\vb{z}}\times\vb{u}  
     +\frac{1}{\fr^2} \grad (\eta - b) - \nu \Delta \vb{u} - \mathbb{F} - \mathbb{B}.$$
     \item The operator $\mathcal{P}: C^1(\Omega, \mathbb{R}) \times C^1(\Omega, \mathbb{R}^2) \to C^0(\Omega, \mathbb{R})$ governs the evolution of the height of the fluid  $\eta$
     $$\mathcal{P}(\eta, \vb{u}):=  \div{(\eta\vb{u})}.$$
    \item $f \in \mathbb{R}$ is the Coriolis parameter, $f=2\Theta \sin\varphi$ where $\Theta$
    is the rotation rate of the Earth and $\varphi$ is the latitude; $f\hat{z}
    \times \vb{u} = (-fv, fu)^T$, where $\hat{z}$ is a unit vector pointing away
    from the centre of the Earth
    \item $\fr = \frac{U}{\sqrt{gH}}$ is the Froude number (dimensionless) which is connected to the stratification of the fluid flow. In this case $U$ is a typical scale for horizontal speed and $H$ is the typical vertical scale, while $g$ is the gravitational acceleration.
    \item $\ro = \frac{U}{f_0 L}$ is the Rossby number (also dimensionless) which describes the effects of rotation on the fluid flow: a small Rossby number ($\ro << 1$) suggests that the rotation term dominates over the advective terms.
    \item $b(x,t)$ is the bottom topography function.
    \item $\nu$ is the viscosity coefficient
    \item The wind forcing is given by $\mathbb{F} = (F, 0)$
with
\[
F(x) = \cos\left( \frac{2\pi x}{L} - \pi \right) + 2 \sin\left( \frac{2\pi x}{L} - \pi \right)
\]
where $x$ is the position in the $x$ direction and $L$ is the domain length in the $x$ direction.
\item The bottom friction $\mathbb{B}$ is represented as a quadratic drag $\mathbb{B}=-\frac{c_D}{\eta}|\mathbf{u}| \mathbf{u}$ with dimensionless drag coefficient $c_D=10^{-3}$. 

\end{itemize}

\begin{figure}[ht!]
\centering
\begin{subfigure}{0.3\textwidth}
    \centering
    \includegraphics[width=\textwidth]{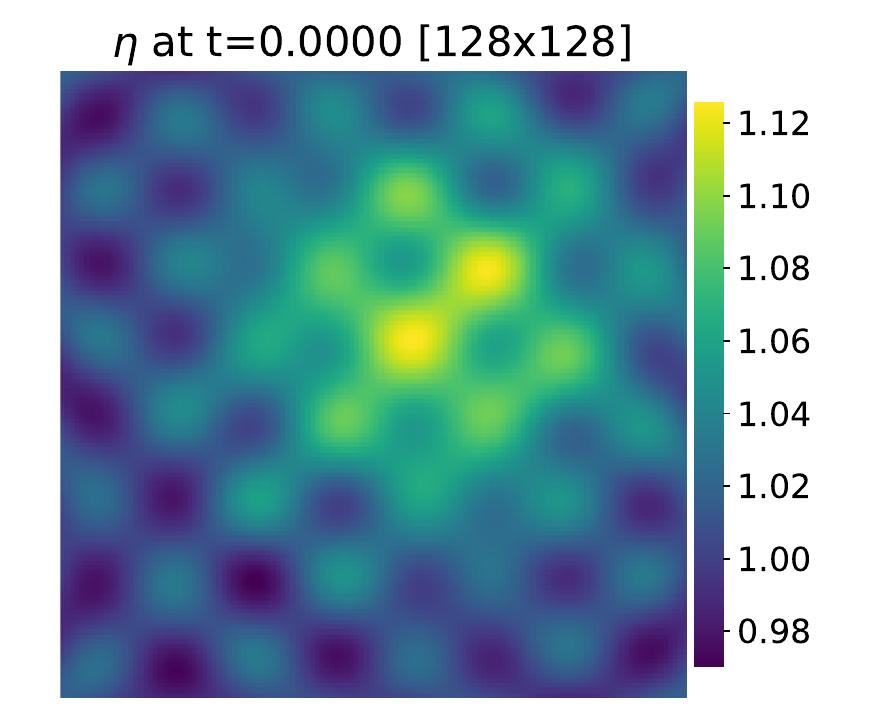}
    \caption{Height~$\eta$.}
    \label{fig:init_fine_e}
\end{subfigure}
\begin{subfigure}{0.3\textwidth}
    \centering
    \includegraphics[width=\textwidth]{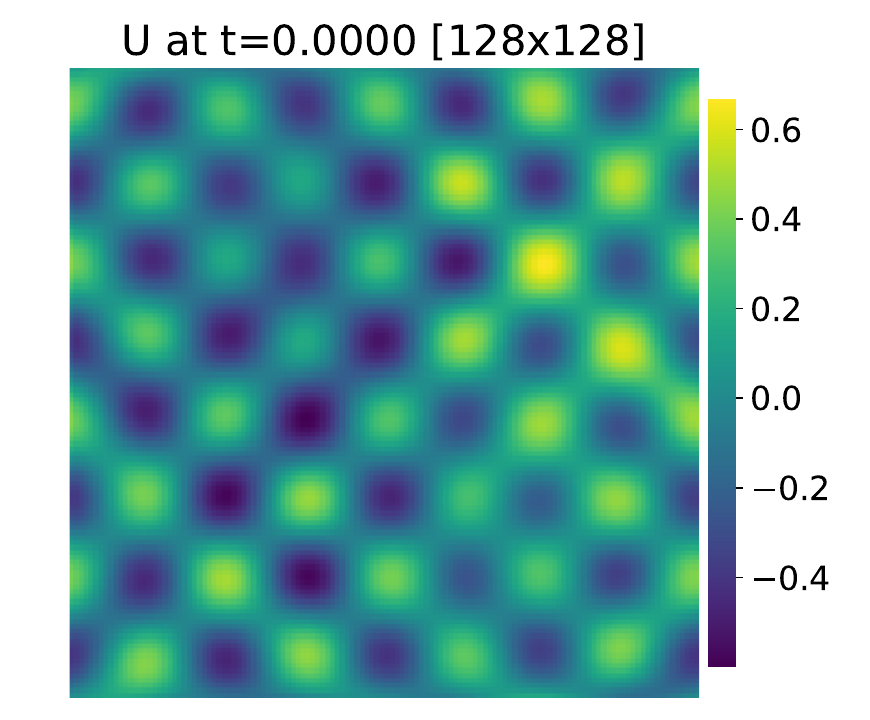}
    \caption{Zonal Velocity~$U$.}
    \label{fig:init_fine_u}
\end{subfigure}
\begin{subfigure}{0.3\textwidth}
    \centering
    \includegraphics[width=\textwidth]{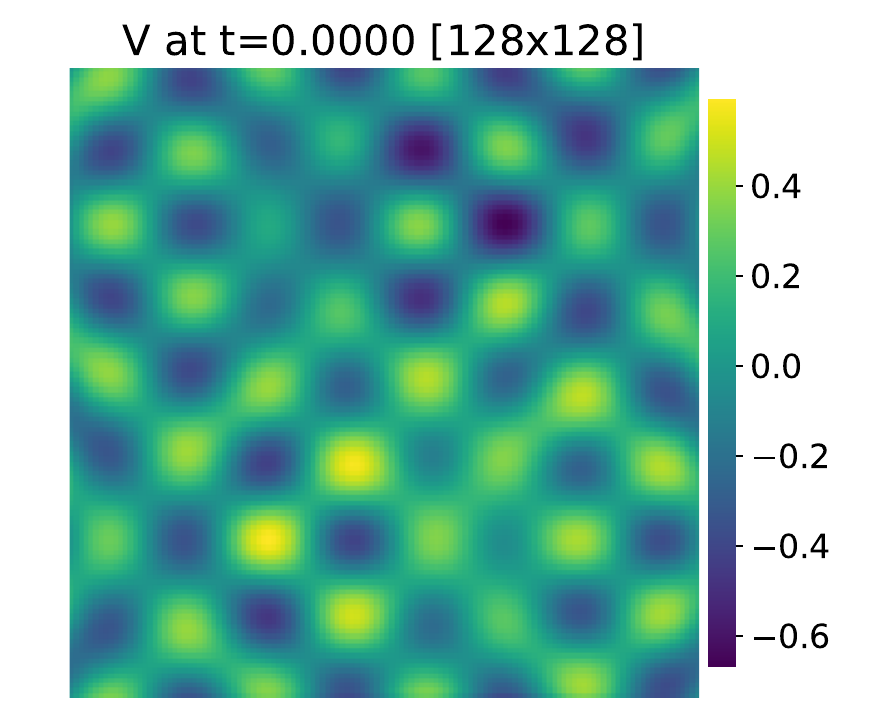}
    \caption{Meridional Velocity~$V$.}
    \label{fig:init_fine_v}
\end{subfigure}
\caption{\textbf{Initial Condition.} Initial condition for the shallow water simulations.} 
\label{fig:rswinit}
\end{figure}

The initial  $\eta$-field is given by  
\begin{equation}
\eta_0 = a \left( 0.05 \, \omega + 0.2 \, e^{-\frac{(x - 0.6)^2 + (y - 0.6)^2}{0.1}}\right) + 1 
\label{eq:init_cond}
\end{equation}
where
\begin{equation}
\omega = \sin(8\pi x) \sin(8\pi y) + 0.4 \cos(6\pi x) \cos(6\pi y) + 0.3 \cos(10\pi x) \cos(4\pi y) + 0.02 \sin(2\pi y) + 0.02 \sin(2\pi x)
\end{equation}
where {$a$ is a free parameter (chosen as $a=0.5$ later on)}. From this initial state corresponding to $\eta$, the starting conditions for the two components of the velocity vector fields are computed using a geostrophic balance condition as
\begin{equation}
    \vb{u}_0 = -\frac{1}{f}\frac{\text{Ro}}{\text{Fr}^2}\grad^\perp \eta_0.
\end{equation}
To achieve $\vb{u}_0 \sim \mathcal{O}(1)$ consistent with the nondimensionalisation, the initial velocity components are globally scaled by the components of $1.5(\max |u|, \max |v|)^T$.
This represents an approximate balance in the system, equating the Coriolis force with the pressure gradient. This geostrophic balance condition is lifted after initial time $0$.
The starting condition is depicted in Figure \ref{fig:rswinit}.
We work with the corresponding discrete version of system (\ref{rsw}), that is
\begin{equation}\label{rswd11}
    \begin{aligned}
     \vb{u}_{n+1}-\vb{u}_{n} + \mathcal{R}(\vb{u}_n, \eta_n)\delta t& = 0\\
     \eta_{n+1}-\eta_{n} + \mathcal{P}(\eta_n, \vb{u}_n)\delta t  &=0.
    \end{aligned}
\end{equation}
In \cite{us} we perturbed the iteration corresponding to (\ref{rswd11}) with spatial Gaussian noise of the form 
\begin{equation}\label{noise}
W_n(x)=\sqrt{\delta t} \sum_{i=1}^{M} \boldsymbol{\xi}_i(x) W_n^i
    \end{equation}
where $(\boldsymbol{\xi}_i)_i$ are divergence-free elements of the covariance basis functions of the SALT (\cite{Holm2015}) noise parametrisation and  $W^i_n\sim N(0,1)$ are independent i.i.d. random variables. When we do this, we obtain the following recurrence formula

\begin{equation}\label{rswd}
    \begin{aligned}
     \vb{u}_{n+1}-\vb{u}_{n} + (\tilde{\vb{u}}_n\cdot\grad) \vb{u}_n + \vb{u}_n \cdot \grad W_n(x) + \frac{f}{\ro} \hat{\vb{z}}\times\tilde{\vb{u}}_n  
     +\frac{1}{\fr^2} \grad (\eta_n - b)\delta t &= \nu \Delta \vb{u}_n \delta t + \mathbb{F} \delta t + \mathbb{B} \delta t\\
     \eta_{n+1}-\eta_{n} + \div{(\eta_n\tilde{\vb{u}}_n)}  &=0
    \end{aligned}
\end{equation}
 where \footnote{The term ${\tilde{\vb{u}}}$ is a velocity perturbation which is specific for this stochastic version of the RSW model.}
\begin{equation}\label{pert}
\tilde{\vb{u}}_n=\vb{u}_n \delta t + W_n(x). 
\end{equation}
The choice of the perturbation (\ref{pert}) is such that the iteration (\ref{rswd}) is an approximation of the stochastic partial differential equation 
\begin{equation}\label{eq:increments1}
    \begin{aligned}
     \id\vb{u} + \mathcal{R}(\vb{u},\eta )\id t +
    + \sum_i \left[(\boldsymbol{\xi}_i \cdot \grad) \vb{u} + \grad\boldsymbol{\xi}_i\cdot\vb{u} + \frac{f}{\ro} \hat{\vb{z}}\times\boldsymbol{\xi}_i\right]\sdiff W_t^i 
    &= 0\\
     \id\eta + \mathcal{P}(\eta, \vb{u}) \id t + \sum_i \div{(\eta\boldsymbol{\xi}_i)}\sdiff W_t^i &=0
    \end{aligned}
\end{equation}
  where $\circ$ denotes Stratonovich integration and $W^i$ are standard i.i.d. Brownian motions as before. The Stratonovich stochastic term generates a second order correction when writing the system in It\^{o} form, but this is dealt with using the intrinsic properties of the numerical scheme. We have explained this part in detail in the Appendix of \cite{us}. The choice of noise in \eqref{eq:increments1} is not arbitrary, but it is inspired by the variational principle approach introduced in \cite{Holm2015}.

 To align the rotating shallow water example with the general notation introduced earlier, note that the state \(m^f\) is represented by the pair \((\vb{u}, \eta)\), which satisfies the partial differential system in equation (\ref{rsw}). Since we are concerned only with discrete approximations, we will directly associate \(m^f\) with the solution of equation (\ref{rswd11}). Similarly, \(m^c\) corresponds to the solution of equation (\ref{rswd}). In other words:
\[
m_{t_n}^c:=\left( 
\begin{array}{c}
\vb{u}_{t_n}^c \\ 
\eta_{t_n}^c
\end{array}
\right)
\]
where $(\vb{u}_{t_n}^c, \eta_{t_n}^c)$ solves \eqref{rswd}. 
Then
\begin{equation}\label{opm}
\mathcal{M}(m_{t_n}^c)(\zeta) = \mathcal{M}\left( 
\begin{array}{c}
\vb{u}_{t_n}^c \\ 
\eta_{t_n}
\end{array}
\right)(\zeta) :=\left( 
\begin{array}{c}
\grad \vb{u}_{t_n}^c \cdot \zeta + u_{t_n}^c \cdot \grad \zeta + \frac{f}{\ro} \hat{\vb{z}}\times\zeta\\ 
\grad \eta_{t_n}^c \cdot \zeta.
\end{array}
\right)
\end{equation}
with $\zeta$ typically corresponding to $\xi_i$, \footnote{Here we can observe one more time the challenges posed by transport noise in general (and SALT noise in this particular case) as this always involves calculating derivatives corresponding to both the model variable $m$ and the (noise) variable $\zeta$. In other words, the operator $\mathcal{M}$ and the variable $\zeta$ are inherently intertwined.} 
and
\begin{equation}
\mathcal{A}(m_{t_n}^c) = \mathcal{A}\left( 
\begin{array}{c}
\vb{u}_{t_n}^c \\ 
\eta_{t_n}^c
\end{array}
\right) =\left( 
\begin{array}{c}
\mathcal{R}(\vb{u}_{t_n}^c, \eta_{t_n}^c) \\ 
\mathcal{P}(\eta_{t_n}^c, \vb{u}_{t_n}^c)
\end{array}
\right).
\end{equation}
Based on \cite{us} we have
\begin{equation}\label{eq:increments} 
\hat{m}_{t_{n+1}} - \hat{m}_{t_n}
 \approx \displaystyle\sum_{i=1}^M\mathcal{M}(m_{t_n}^c)\xi_i(x)\sqrt{\Delta} W_{t_n}^i = \mathcal{M}(m_{t_n}^c) W_{t_n}(x)
\end{equation}
with $\mathcal{M}$ given above in \eqref{opm} and
\begin{equation}
W_{t_n}(x)=\sqrt{\delta t} \sum_{i=1}^{M} \boldsymbol{\xi}_i(x) W_{t_n}^i
\end{equation}
for the particular case of the RSW model.
 In previous work \cite{us}, we first generated the increments in \eqref{eq:increments} for $\eta$ only and then we used them together with a geostrophic balance assumption to compute the corresponding noise increments for the two components of $\vb{u}_{t_n}^c$. This means that for the RSW model we work mainly with $\mathcal{M}(\eta_{t_n}^c)(\zeta) = (\grad \eta_{t_n}^c \cdot \zeta)$ which corresponds to a \textit{transport noise}. The approach we adopted in \cite{us} was based on a Principal Component Analysis (PCA) methodology that was undepinned by the ansatz that the increments are normally distributed. This ansatz is however not justified in the current case, as one can see in Figure \ref{fig:distribution_histograms} (a) where the distribution of the noise is bimodal.  
In \cite{us2} we replaced the spatial Gaussian noise $W_{t_n}(x)$ in \eqref{eq:increments} with a general noise $N_n(x)$ 
such that
\begin{equation}\label{diffnew}
\hat{m}_{t_{n+1}}-\hat{m}_{t_{n}}\approx \displaystyle\mathcal{M}%
(m_{t_{n}}^{c})N_{n}  
\end{equation}
{where }$N_{n}$ has an unknown distribution which is modelled using a diffusion Schr\"{o}dinger bridge.
In the current work, we proceed with the data assimilation procedure using a signal which is calibrated based on the generative model methodology described in \cite{us2}.

\section{Methodology}

In the following we give a brief description of the calibration approach, which is based on the Schrödinger bridge, and the particle filtering methodology for performing data assimilation.

\subsection{Noise Calibration }\label{sect:calib}
\begin{figure}[ht!]

\centering
\begin{subfigure}{0.3\textwidth}
    \centering
    \includegraphics[width=\textwidth]{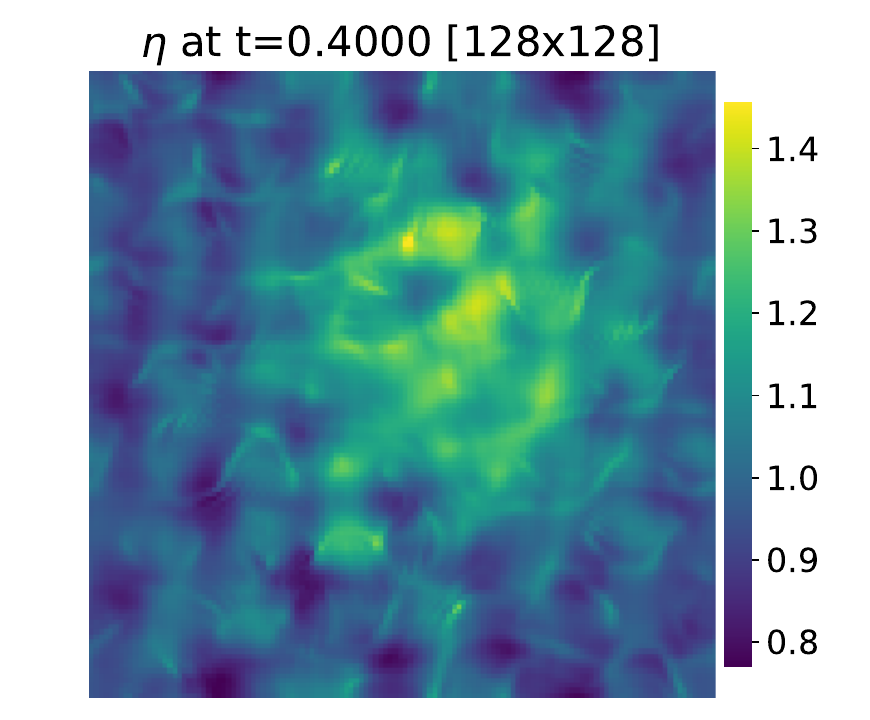}
    \caption{Height $\eta$ Fine.}
    \label{fig:var_fine_e}
\end{subfigure}
\begin{subfigure}{0.3\textwidth}
    \centering
    \includegraphics[width=\textwidth]{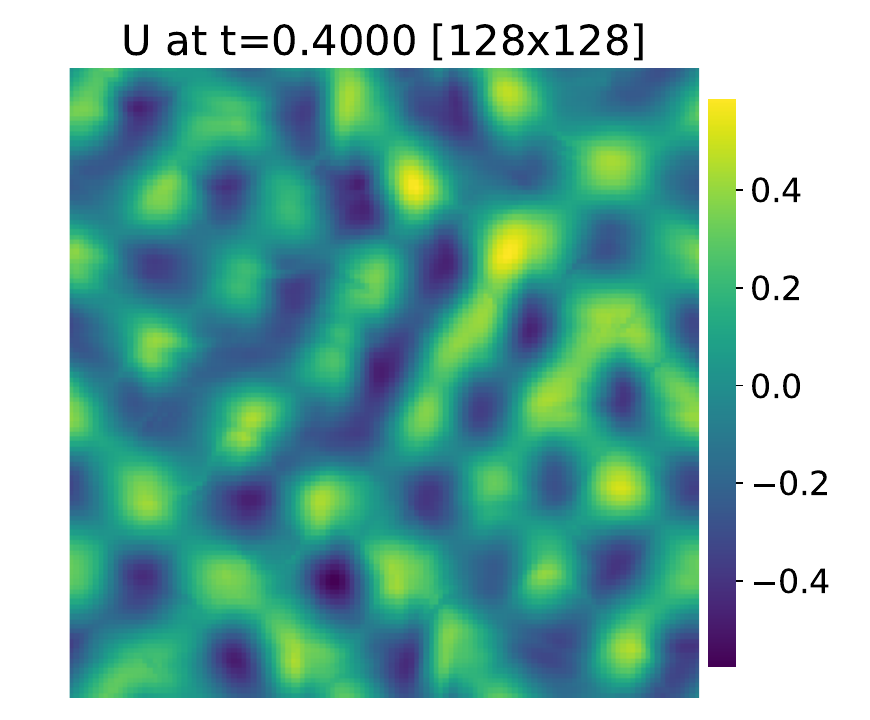}
    \caption{Zonal Velocity~$U$ Fine.}
    \label{fig:var_fine_u}
\end{subfigure}
\begin{subfigure}{0.3\textwidth}
    \centering
    \includegraphics[width=\textwidth]{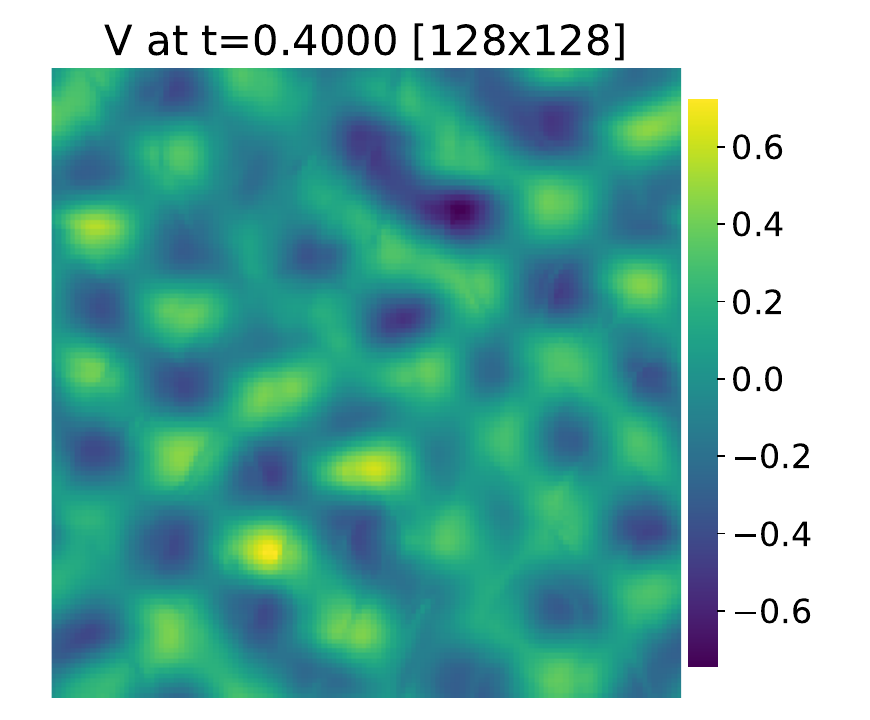}
    \caption{Meridional Velocity~$V$ Fine.}
    \label{fig:var_fine_v}
\end{subfigure}
\hfill 
\begin{subfigure}{0.3\textwidth}
    \centering
    \includegraphics[width=\textwidth]{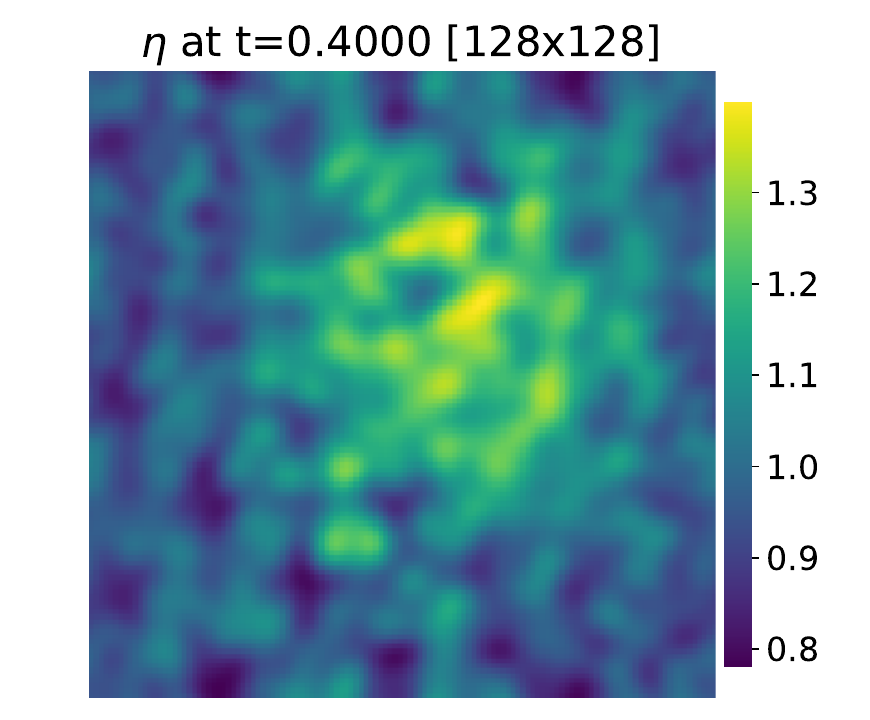}
    \caption{Height $\eta$ Coarse.}
    \label{fig:var_coarse_e}
\end{subfigure}
\begin{subfigure}{0.3\textwidth}
    \centering
    \includegraphics[width=\textwidth]{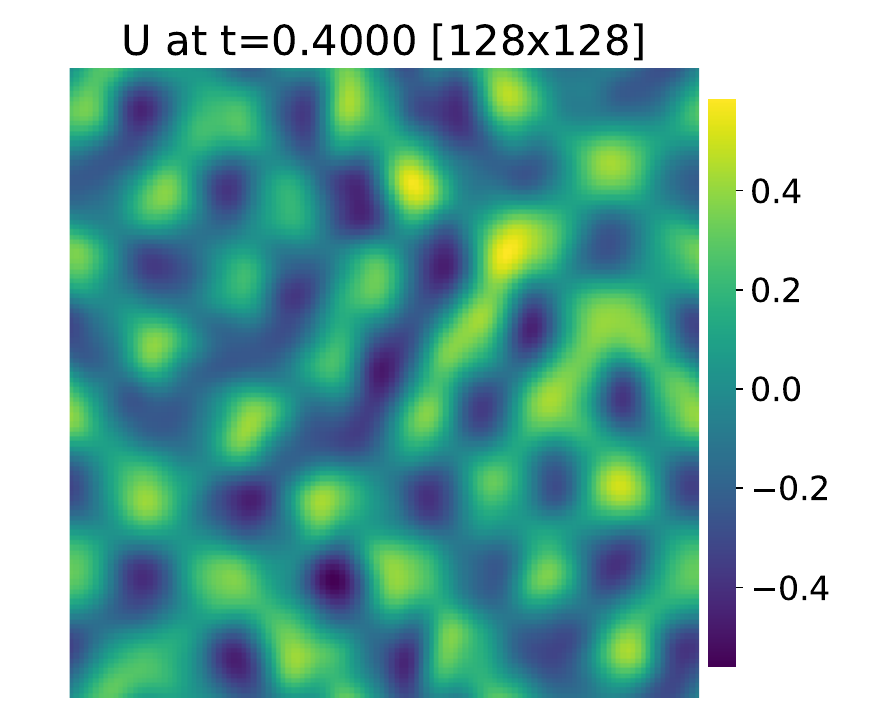}
    \caption{Zonal Velocity~$U$ Coarse.}
    \label{fig:var_coarse_u}
\end{subfigure}
\begin{subfigure}{0.3\textwidth}
    \centering
    \includegraphics[width=\textwidth]{figures/var_U_grid_128_T_4000_coarsened.pdf}
    \caption{Meridional Velocity~$V$ Coarse.}
    \label{fig:var_coarse_v}
\end{subfigure}
\caption{\textbf{Shallow Water Snapshots.} Plots of the full and filtered (rows) non-dimensional shallow water variables (columns). The variables (a+d) height~$\eta$, (b+e) zonal velocity~$U$, and (c+f) meridional velocity~$V$ are plotted on the fine grid of resolution $128\times 128$ grid points. The top row (a-c) shows the fully resolved fields and the bottom row (d-f) shows the coarsened fields after applying the low-pass filter with a maximum frequency of $16$ gridpoints. The solution is obtained after running $4000$ timesteps up to time $t=0.4$ started from the prescribed initial condition.} 
\label{fig:PDE_plot}

\end{figure}

We are calibrating the Nondimensional Rotating Shallow Water System on the unit square $\Omega=[0,1]\times[0,1]$. In Figures \ref{fig:PDE_plot} (a)-(c) we see the state of the system after 4000 time steps for the three variables corresponding to the rotating shallow water model, on a fine resolution grid of size 128 x 128. 

 \paragraph{Mollification of the fine grid data.}
{{We use the $\eta$-component data of the PDE run on the fine grid to obtain $\eta^f$ and construct a coarsening $C(\eta^f)$ by mollifying $\eta^f$.}}
For the mollification, we chose to apply a low-pass filter. Consider the field $\eta^f$ obtained as above.
Its 2D Fourier transform in space is denoted by $\hat{\eta}^f$. The mollified field is then obtained by truncating $\hat{\eta}^f$ such that 
\begin{equation*}
    C(\eta^f) = \operatorname{IFT} \circ \iota \circ \hat{\eta}^f
\end{equation*}
where $\operatorname{IFT}$ denotes inverse Fourier transform and  $\iota$ is the function such that
\begin{equation*}
    \iota( r)(\omega) = 
    \begin{cases}
        r(\omega) &\text{if}\;\; \omega < f_{max}, \\
        0 &\text{if}\;\; \omega \geq f_{max}.
    \end{cases}
\end{equation*}
Here, $f_{max}$ denotes a cutoff frequency. Any spatial frequencies above $f_{max}$ present in $\eta^f$ are removed. 

As mentioned above, this procedure is known as a low-pass filter in signal processing. In our context, the low-pass filter is a principled choice, since the high frequency effects present in the field will be removed which corresponds to the fact that we expect coarser simulations to misrepresent higher frequency effects. Alternative coarsenings may be performed using Gaussian mollifiers, spatial averaging, or subsampling.

In Figures \ref{fig:PDE_plot} (d) - (f) we see the mollification corresponding to the rotating shallow water system after applying a low-pass filter with a maximum frequency of 16 grid points and running for 4000 steps.  
  
\paragraph{Data generation.}  We prepare the data required to estimate 
the noise parametrisation.  We use the time-increments of the discrepancy between the high resolution fields and coarsened fields. This corresponds to an Euler-Maruyama step of the stochastic integral. 
The samples from the noise increments are computed from this data by solving a sequence of hyperbolic equations\footnote{Note that these calibration equations do not have a unique solution. We have not explored fully the effect that this non-uniqueness result on the calibration results.}  
 \begin{equation}\label{calibc}
 \delta \hat{\eta_i} = C(\eta_{i,x}^f)\frac{\partial \psi_i}{\partial y}-C(\eta_{i,y}^f)\frac{\partial \psi_i}{\partial x}
 \end{equation}
for each calibration time/datapoint indexed by $i$, and with periodic boundary condition. The subscript $x$ and $y$ denote partial differentiation.
 The solutions $\psi_i$ to the sequence of calibration equations can be thought of as \emph{stream functions} for the perturbation fields $\tilde{u}$, $\tilde{v}$, so that $(\tilde{u},\tilde{v})^\top = \nabla^\perp \psi$. The 2D grid based data $(\tilde{u}_{ij},\tilde{v}_{ij})$ is vectorized as $\Psi = ((\tilde{u}_{ij}), (\tilde{v}_{ij}))^\top$ and represented in the following form 
\begin{equation*}
{\Psi _{i}-\bar{\Psi}}=\Delta N_i,
\end{equation*}%
where $\Delta N_i$ are i.i.d. random variables that are the training data for the generative model. A plot of the distribution of the training data is shown in Figure~\ref{fig:data_hist}

\subsection{Diffusion Schrödinger Bridge}

In this work, we use the Diffusion Schrödinger Bridge (DSB) model~\cite{de2021diffusion}, a particular type of score-based generative models, to calibrate the hybrid shallow water system introduced above. Score-based generative models are considered state-of-the-art in generative modeling and have demonstrated strong performance across various benchmark tasks in machine learning \cite{vahdat2021score, song2020improved, zimmermann2021score}. We employ the learnt model to generate new samples from the underlying probability distribution of synthetic noise increments which are derived from the calibration equations \eqref{calibc}. These constitute the training set data. Neural network-based generative models are typically chosen for this task due to the high-dimensional nature of the data distribution.

The fundamental principle of diffusion generative models is that data is gradually diffused by incrementally adding noise, until it resembles a sample from a pure noise distribution (e.g. a Gaussian distribution). A backward process which is rigorously defined is then used to generate samples from the unknown data distribution. This involves drawing a sample from the pure noise distribution and then executing the reverse diffusion process, known as the \emph{denoising} process. 

In general, the DSB model introduced in~\cite{de2021diffusion} enhances the classical score-based diffusion model by incorporating an optimal transport procedure called \emph{iterative proportional fitting (IPF)}. This allows for iterative score-based diffusion training and enables the use of shorter noising and denoising processes. Additional details regarding the application of the DSB model in the calibration procedure can be found in \cite{us}.

\subsection{Particle Filters}\label{sect:pf}

The standard particle filter performs effectively for models of small to medium size. However, as the dimensionality of the model increases, the filter's performance deteriorates. Intuitively, with higher dimensions, the gap between the predictive and posterior distributions widens. Consequently, particles may deviate significantly from the true state (due to the larger space in higher dimensions) and could end up having very low likelihoods, possibly decreasing exponentially. To address this issue, a new methodology has been introduced in \cite{CrisanBeskosJasra}. Instead of a single resampling step, the standard particle filter undergoes repeated cycles of tempering, resampling, and jittering, with the number of cycles being adaptively determined based on the effective sample size of the intermediate sample. As demonstrated in \cite{CrisanBeskosJasra}, this approach stabilizes the particle filter in high dimensions, preventing its degeneration even as the dimensionality grows, while maintaining a fixed number of particles. In this paper, we successfully apply this methodology to the stochastic hybrid rotating shallow water model.
We give next a brief description of each of the three steps:

\emph{Tempering} introduces a series of intermediate probability distributions that allow for a smooth transition between the prior (predictive) distribution and the posterior (target) distribution. As a result, the filter can gradually incorporate new observational data rather than doing it all at once, which helps to prevent particles from drifting too far away from the truth. Therefore, instead of directly updating the particles using the full likelihood of the observations, tempering defines a sequence of intermediate distributions, with a parameter that controls the interpolation between prior and posterior. At each intermediate step, the particles are weighted according to the likelihood relative to this intermediate distribution, helping them stay closer to the true posterior as more data is assimilated.

\emph{Resampling} addresses particle degeneracy, where only a few particles with large weights remain after several iterations, leading to reduced diversity and performance. Following tempering, particles are resampled according to their weights: higher-weight particles are duplicated, while those with low weights are discarded. This process maintains a particle population that better represents the posterior distribution, ensuring effective filtering.

\emph{Jittering} adds an MCMC perturbation to the particles after resampling to increase the particle diversity. This prevents the filter from collapsing into a few points and keeps the particle population more spread out, maintaining a better exploration of the state space. After resampling, some noise is introduced into the resampled particles to "jitter" or perturb their positions slightly. This ensures that duplicated particles from the resampling step are not identical and can continue to explore the state space effectively.

In short, the tempering step helps the filter transition smoothly between the prior and posterior, the resampling step keeps the particles relevant, and the jittering step maintains diversity in the particle population. This sequence of steps helps ensure that the particle filter remains stable and effective even when the dimension of the model is large.

\section{Numerical Results}
In the following we present the numerical results. The results of the calibration experiments are given in Subsection~\ref{sec:num_calib}. The results for the data assimilation experiments are explained in Subsection~\ref{sec:num_da}. All simulations of the non-dimensional shallow water model are run using a Froude number of $\fr = 1.1$ and a Rossby number of $\ro = 0.2$.
\subsection{Calibration}\label{sec:num_calib}

Data for the calibration of the SPDE noise was generated from a simulation of the deterministic model on a fine grid of size $128 \times 128$ grid points over a total duration of $15,000$ timesteps up to a final time of $t=1.5$.
The data is collected by applying the low pass filter with a maximum frequecy $k_{\text{max}}=16$ corresponding to a coarse grid resolution of $32\times 32$ gridpoints.
We isolate the high-frequency fluctuations by taking the difference between the original and filtered fields. The computed fluctuations of the height variable $\eta$ are then time-differenced to yield increments and used as the right hand side for the data-dependent calibration equation.
A solution of the calibration equation is approximated using a finite-difference method constrained to yield a mean-zero solution.
The approximated solutions represent a stream function for the noise we seek to generate. Iterating this process over time generates a sample from the distribution of the assumed stream function. 

\begin{figure}
\centering
\begin{subfigure}{0.4\linewidth}
    \includegraphics[width=\linewidth]{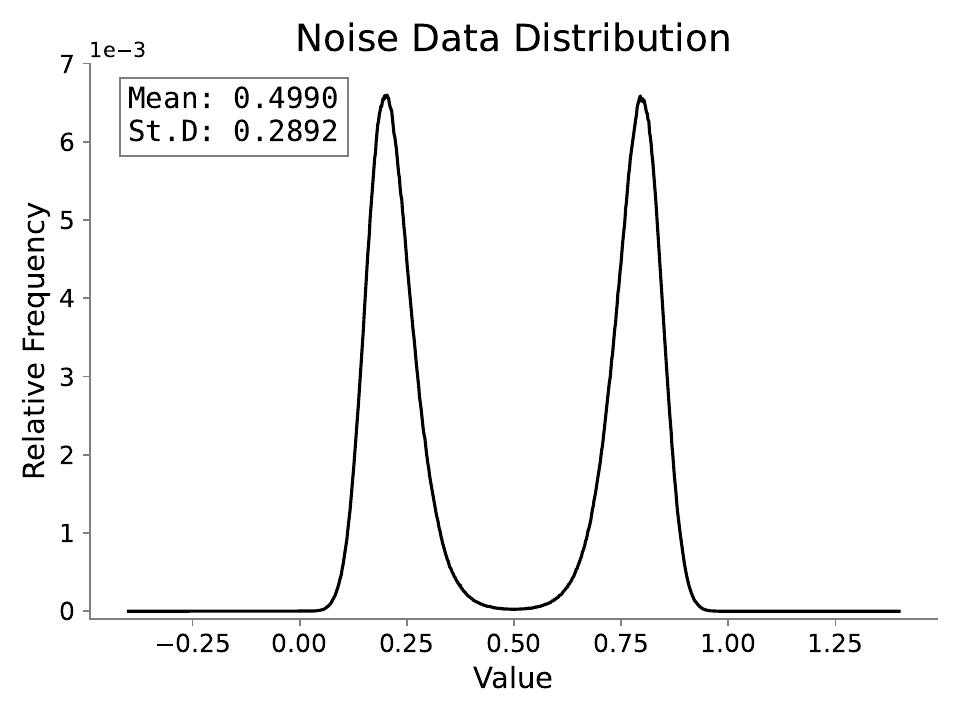}
    \caption{Distribution of the Noise Data}
    \label{fig:data_hist}
\end{subfigure}
\begin{subfigure}{0.4\linewidth}
    \includegraphics[width=\linewidth]{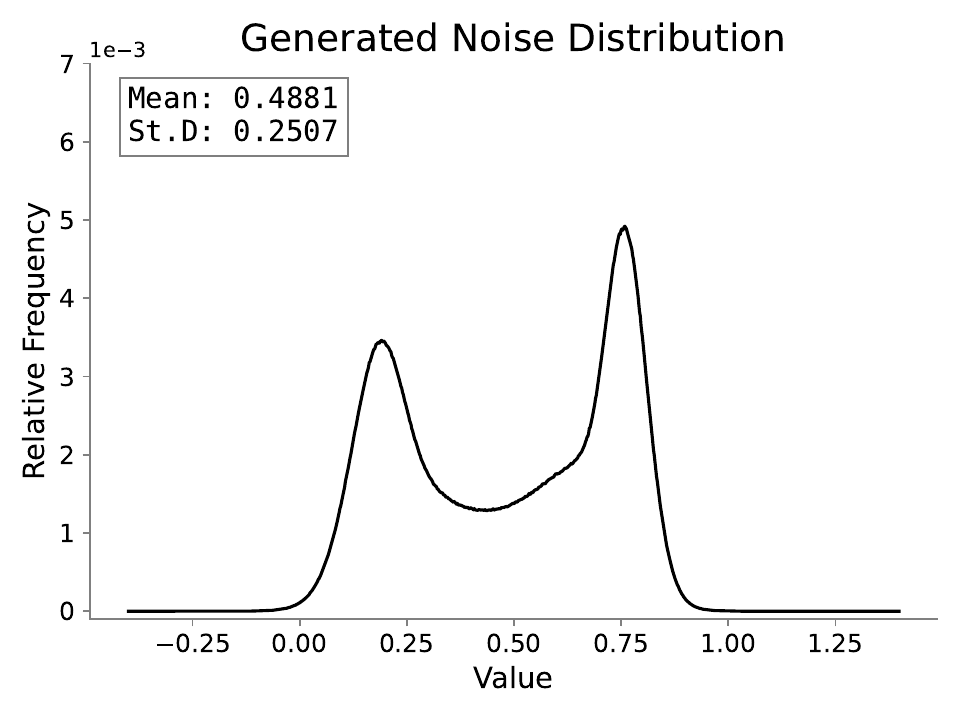}
    \caption{Distribution of the Generated Noise}
    \label{fig:gen_model_hist}
\end{subfigure}
\caption{Distribution of the Noise Data obtained from (a) the deterministic shallow water model and produced by (b) the trained generative model.}
\label{fig:distribution_histograms}
\end{figure}
A reversible statistical transformation is then applied to this dataset, by first applying an $\operatorname{arcsinh}$-transform~\cite{burbidge1988alternative} and subsequently normalising the data to the range $[0,1]$. A plot of the distribution of the dataset across all samples and all grid locations is depicted in Figure~\ref{fig:data_hist}. The distribution is clearly bimodal and approximately symmetric around the mean.
Then the generative model is trained on this dataset. The training dataset consists of $N_{\text{train}} = 14,999$ samples. The model is trained using $30$ timesteps in the diffusion process and over $9$ DSB steps with $10,000$ training iterations per DSB step.

\begin{figure}
    \centering
    \includegraphics[width=0.5\linewidth]{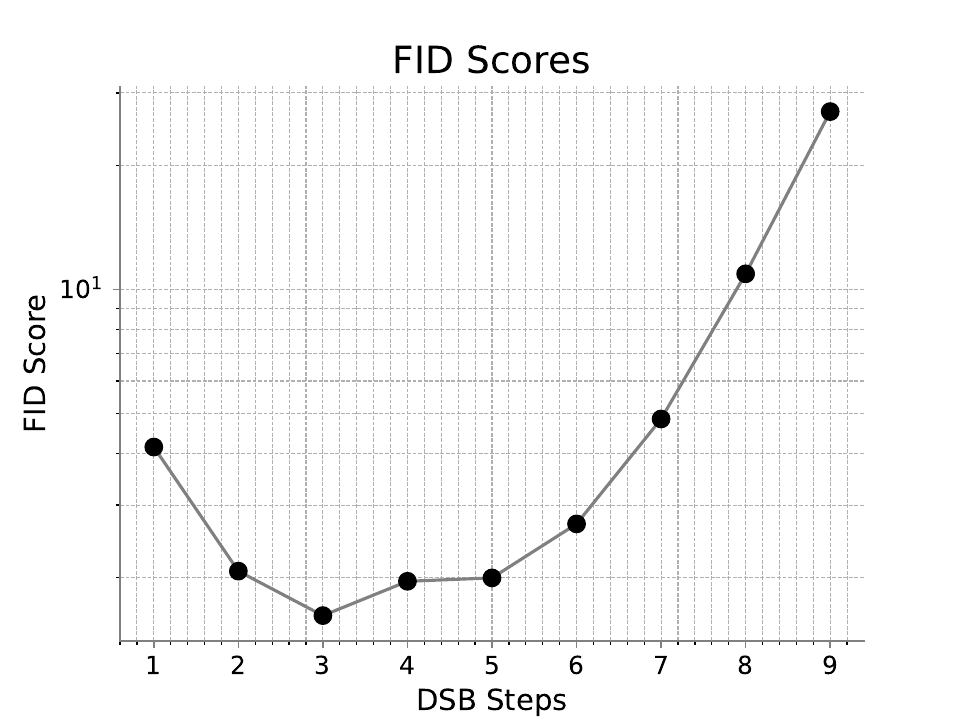}
    \caption{FID Scores over the number of DSB steps}
    \label{fig:fid}
\end{figure}

We identify the best model in terms of the FID score (see \cite{heusel2017gans}) as the one after $3$ DSB steps. Also, the models after $2$, $4$, and $5$ DSB steps achieve very good scores. All subsequent runs use data generated by the DSB 5 model.

After training, we generate a set of $N_{\text{noise}}= 40,000$ noise samples by evaluating the trained generator. The generated noise set serves as a dictionary for the noise used in the following ensemble runs. The distribution across samples and space of the samples generated by the diffusion model is shown in Figure~\ref{fig:gen_model_hist}. We observe that the generative model clearly reproduces the bimodality of the original dataset, placing the modes in approximately the same locations as in Figure~\ref{fig:data_hist}. The distribution of generated samples appears slightly less symmetric than the original dataset, with more weight placed around the mean.

The generated noise is used to run forecast ensembles using the stochastic RSW model on the coarse grid of $32 \times 32$ grid points. To inject the noise into the SPDE it was necessary to clip the range of the noise samples to $\pm 1$ standard deviation around the local mean at each grid point. This is done to avoid large gradients when computing the velocity perturbations from the stream function samples. Additionally, to control the strength of the noise, we introduce a scaling factor for the noise samples set to $30$ throughout the following simulations.

We run the forecast for three different forecast horizons over a total of $4000$ steps. We use a forecast horizon of $1000$ steps repeated $4$ times over the total duration and reset to relaunch from the truth each time. Additionally, we run a horizon of $2000$ timesteps repeated twice, and finally a forecast with horizon $4000$.
In each case the ensemble is started from the truth. Here, the truth is the deterministic model run on the fine $128 \times 128$ grid and evaluated on the coarse grid locations. The ensemble forecasts and observations of the truth are collected at the final times of the forecasts on all $1024$ points in the computational domain. Each of the variables $\eta$, $U$, and $V$ are treated separately.

\begin{table}[ht!]\small
    \centering
    \caption{\textbf{CRPS.} CRPS Scores for all variables (columns) with different forecast horizons (rows). Smaller CRPS scores are better.}
    \begin{tabular}{ c c c c }
    \toprule
        \textbf{Horizon} & $\boldsymbol{\eta}$ & $\mathbf{U}$ & $\mathbf{V}$ \\
        \midrule
        $\mathbf{1000}$ & $0.0614$ & $0.0703$ & $0.0913$ \\
        $\mathbf{2000}$ & $0.0731$ & $0.0596$ & $0.0655$ \\
        $\mathbf{4000}$ & $0.0794 $& $0.0672$ & $0.0757$ \\
        \bottomrule
    \end{tabular}
    
    \label{tab:crps_scores}
\end{table}

We assess the calibration of the noise in the SPDE using various forecast metrics. 
In Table~\ref{tab:crps_scores} we report the continuous ranked probability scores (CRPSs) (see, for example \cite{DecompositionoftheContinuousRankedProbabilityScoreforEnsemblePredictionSystems}) for the different scenarios, averaged over the domain and the forecast repetitions.
We observe that for the variable $\eta$, the CRPS score is best (smallest) for the short horizon of $1000$ steps, whereas the velocity forecasts for $U$ and $V$ are both most effective for the longer horizon of $2000$ timesteps.

\begin{figure}
\centering
\begin{subfigure}{0.25\linewidth}    
    \includegraphics[width=\linewidth]{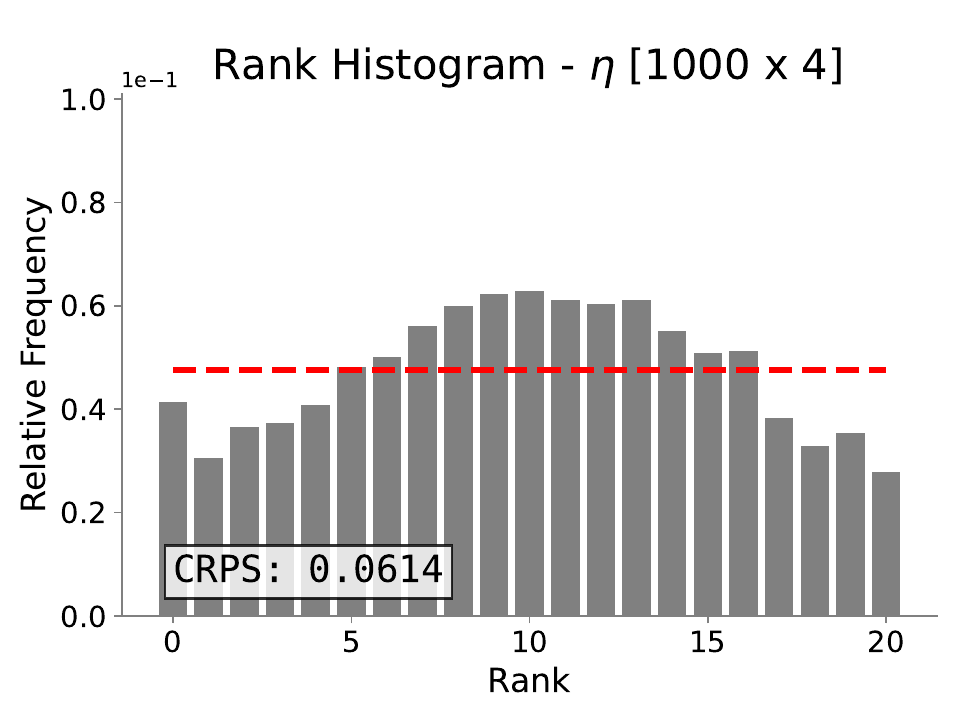}
    \caption{$\eta$ for Horizon 1000}
    \label{fig:rhH1000}
\end{subfigure}
\begin{subfigure}{0.25\linewidth}
    \includegraphics[width=\linewidth]{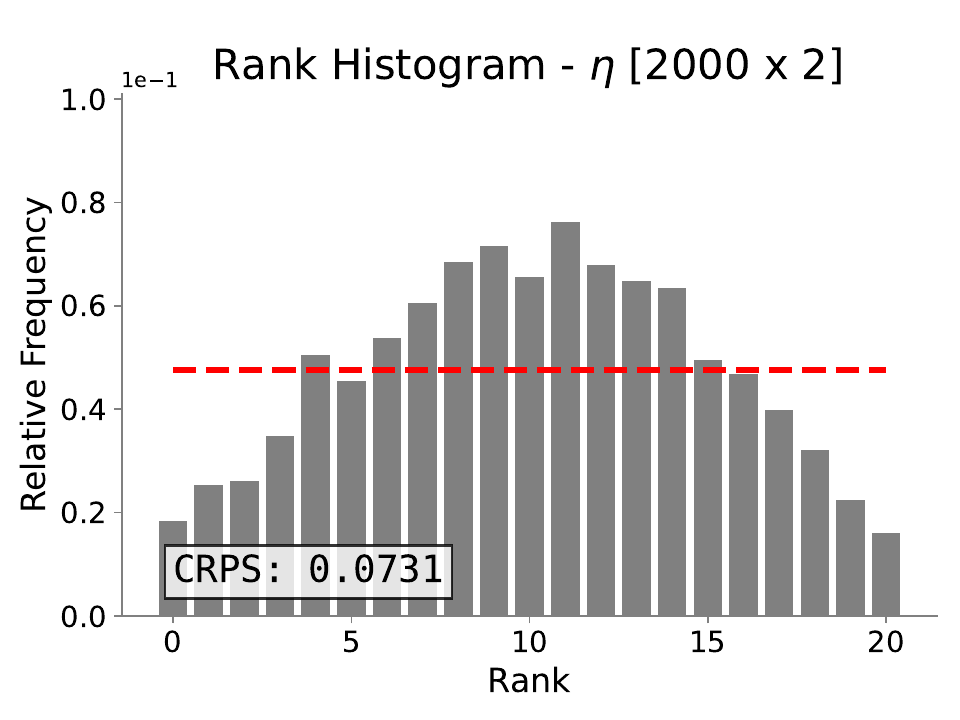}
    \caption{$\eta$ for Horizon 2000}
    \label{fig:rhH2000}
\end{subfigure}
\begin{subfigure}{0.25\linewidth}
    \includegraphics[width=\linewidth]{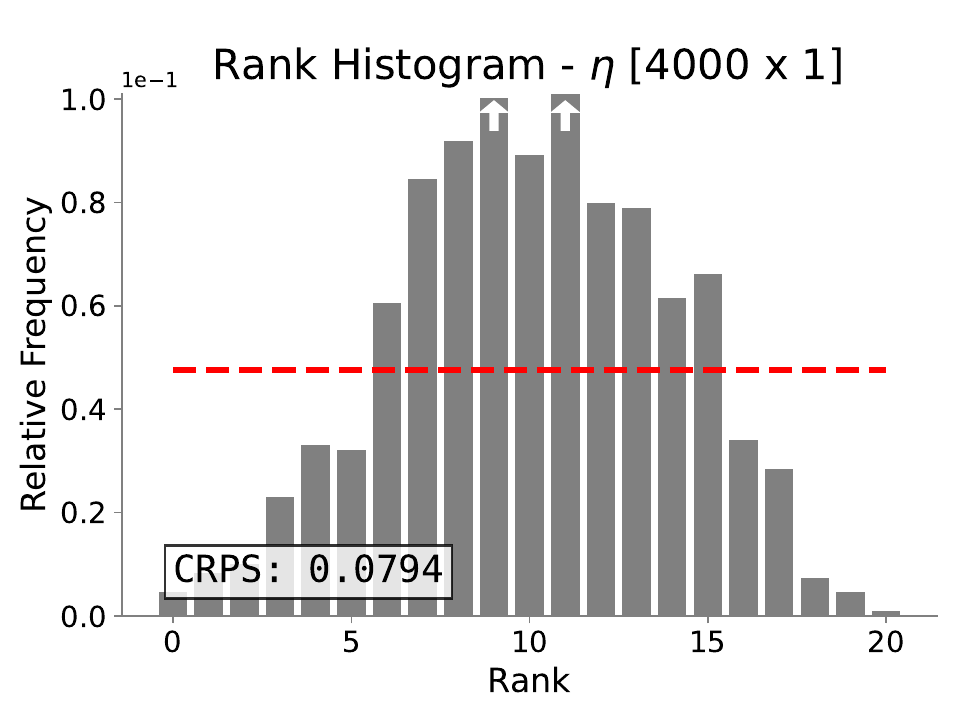}
    \caption{$\eta$ for Horizon 4000}
    \label{fig:rhH4000}
\end{subfigure}
\\
\begin{subfigure}{0.25\linewidth}    
    \includegraphics[width=\linewidth]{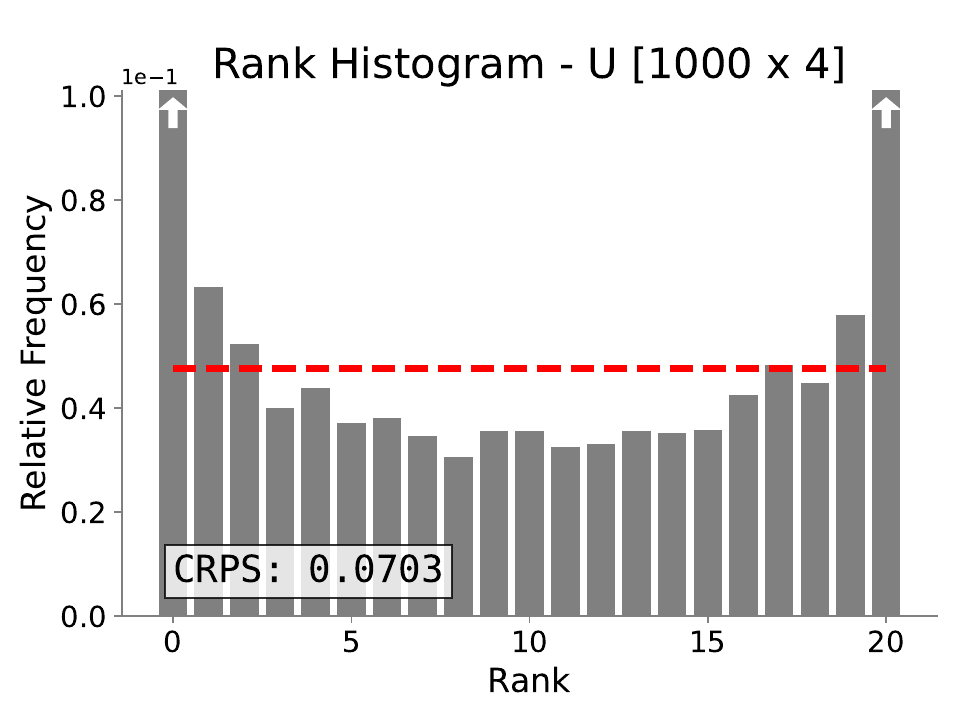}
    \caption{U for Horizon 1000}
    \label{fig:rhU1000}
\end{subfigure}
\begin{subfigure}{0.25\linewidth}
    \includegraphics[width=\linewidth]{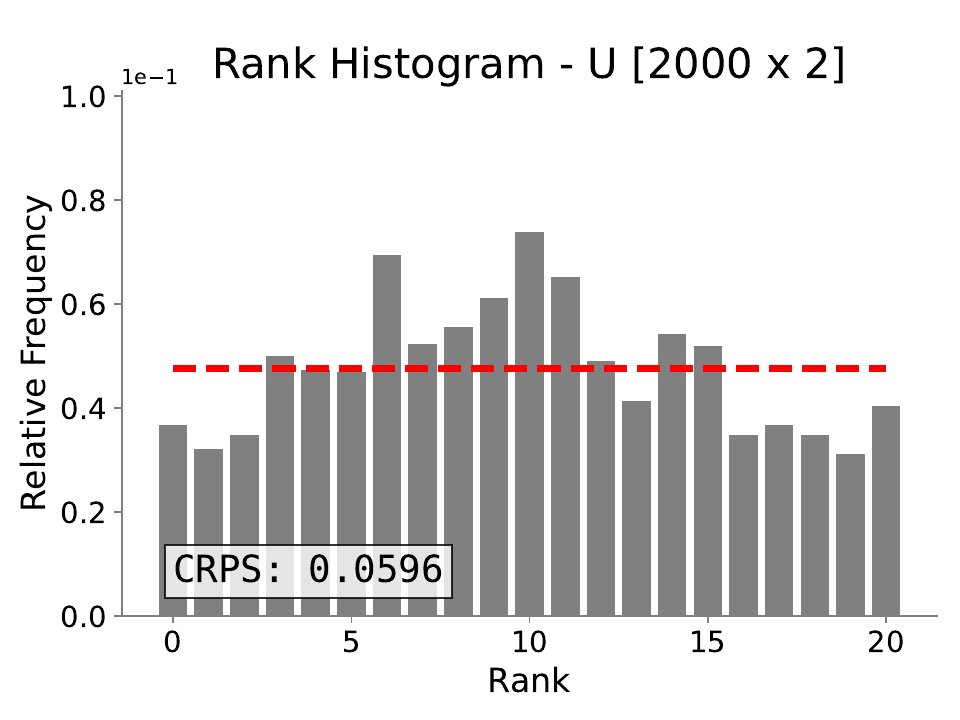}
    \caption{U for Horizon 2000}
    \label{fig:rhU2000}
\end{subfigure}
\begin{subfigure}{0.25\linewidth}
    \includegraphics[width=\linewidth]{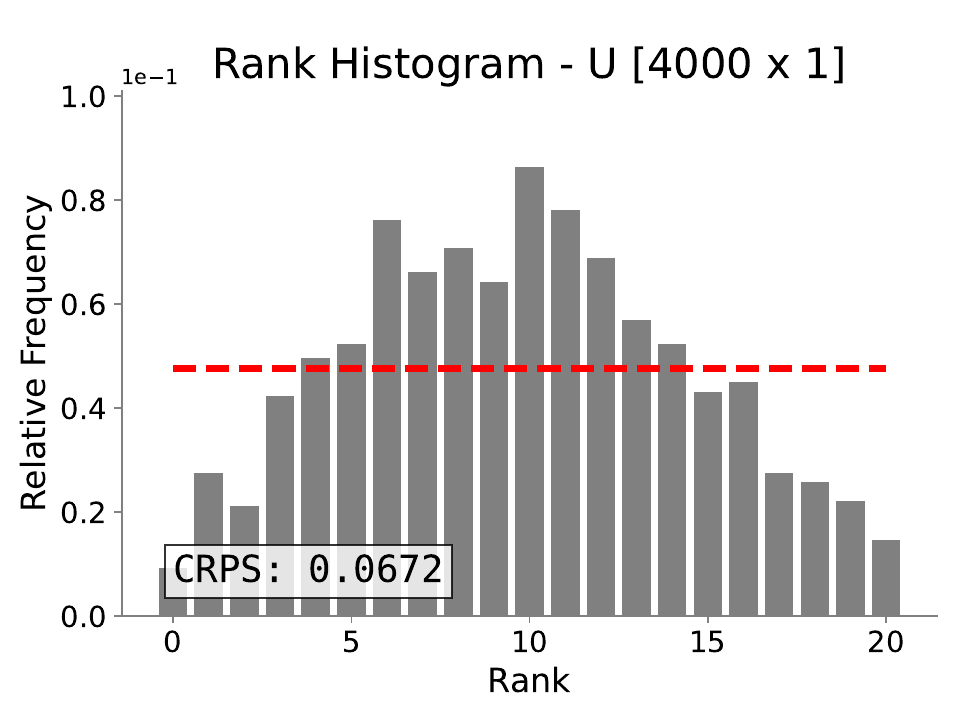}
    \caption{U for Horizon 4000}
    \label{fig:rhU4000}
\end{subfigure}
\\
\begin{subfigure}{0.25\linewidth}    
    \includegraphics[width=\linewidth]{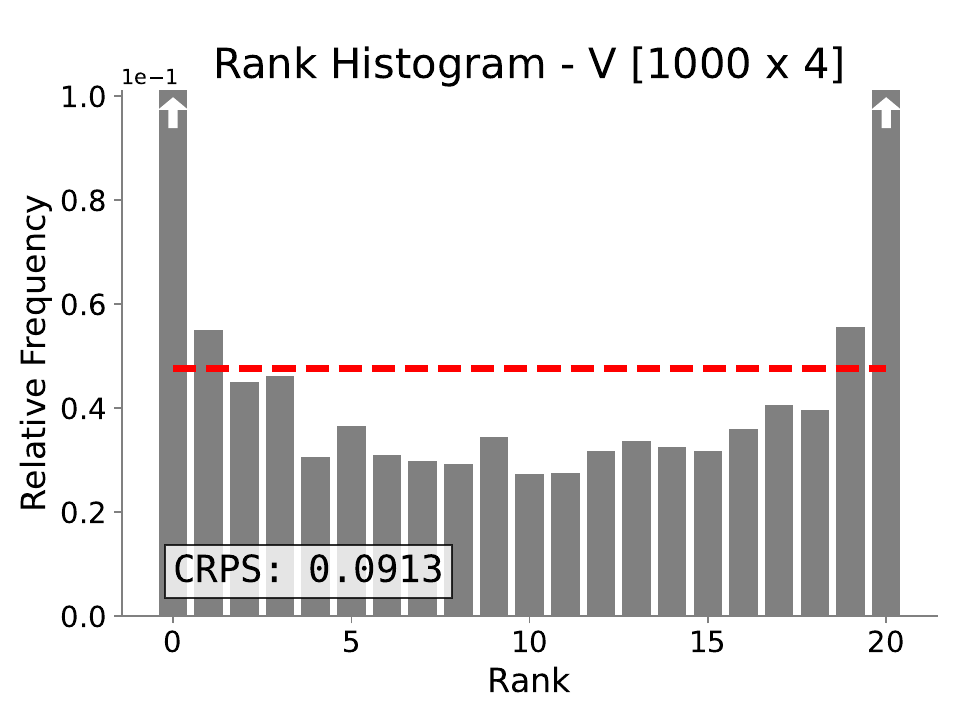}
    \caption{V for Horizon 1000}
    \label{fig:rhV1000}
\end{subfigure}
\begin{subfigure}{0.25\linewidth}
    \includegraphics[width=\linewidth]{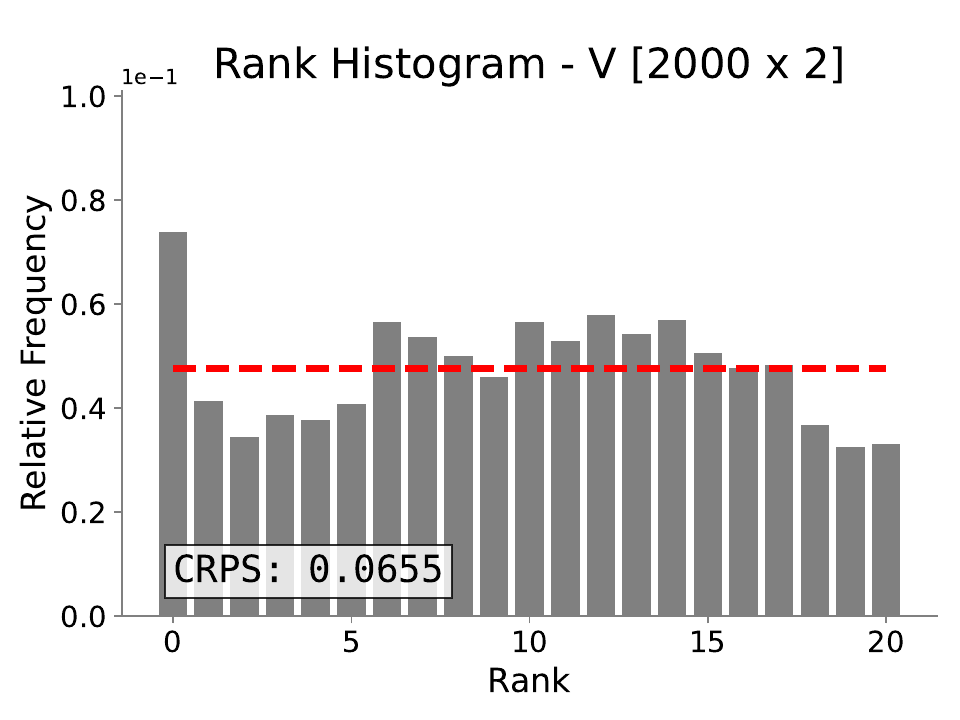}
    \caption{V for Horizon 2000}
    \label{fig:rhV2000}
\end{subfigure}
\begin{subfigure}{0.25\linewidth}
    \includegraphics[width=\linewidth]{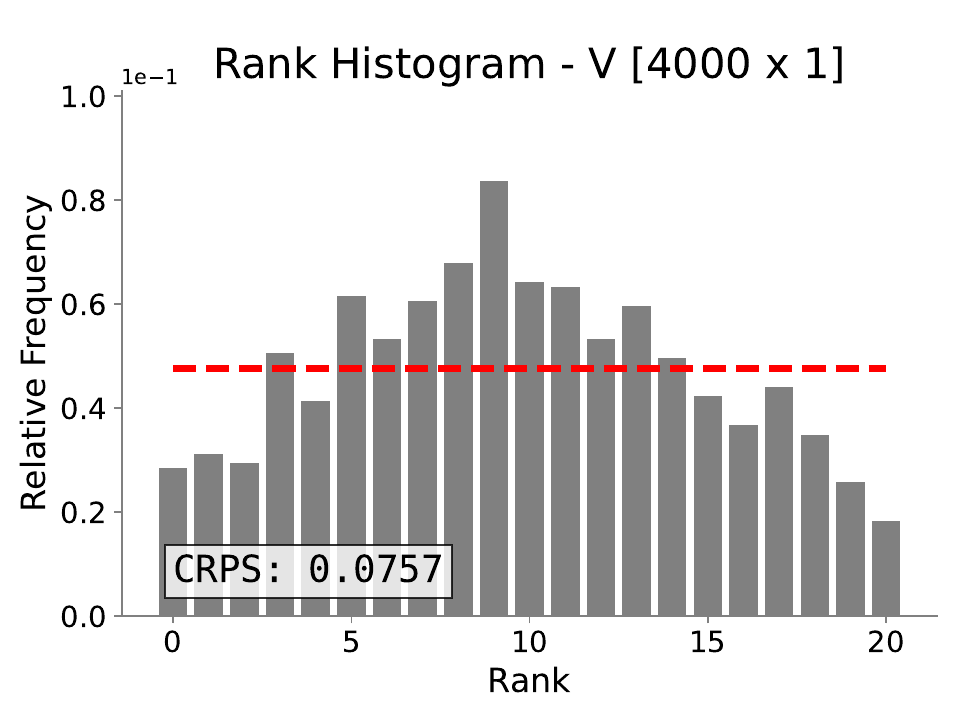}
    \caption{V for Horizon 4000}
    \label{fig:rhV4000}
\end{subfigure}
\caption{{\textbf{Rank Histograms.}} Rank histograms for the variables $\eta$, U and V (rows) and different forecast horizons (columns). Red horizontal lines represent theoretically perfect balance at $1/(N_{\text{ens}}+1)$. White arrows within bars indicate that the relative frequency exceeds $0.1$.}
\label{fig:rank_histograms}
\end{figure}

The next forecast quality study was done in terms of the rank histograms (see, for example \cite{InterpretationofRankHistogramsforVerifyingEnsembleForecasts}) reported in Figure~\ref{fig:rank_histograms}. Visual inspection of the rank histograms confirms the observation from the CRPS scores, namely that the shorter horizon of $1000$ steps yields the most effective forecast for the $\eta$ variable (see Figure~\ref{fig:rhH1000}), whereas the forecast for velocity $U$ and $V$ is most effective for the medium forecast of $2000$ timesteps (see Figures~\ref{fig:rhU2000} and \ref{fig:rhV2000}). We also observe that the forecast ensembles for $\eta$ become gradually more overdispersed for longer horizons (see Figures~\ref{fig:rhH2000} and \ref{fig:rhH4000}). Moreover, the short forecasts for $U$ and $V$ are underdispersed (see Figures~\ref{fig:rhU1000} and \ref{fig:rhV1000}), whereas the long velocity forecasts are each overdispersed (see Figures~\ref{fig:rhU4000} and \ref{fig:rhV4000}).

\begin{figure}
    \centering
    \begin{subfigure}{0.3\linewidth}    
        \includegraphics[width=\linewidth]{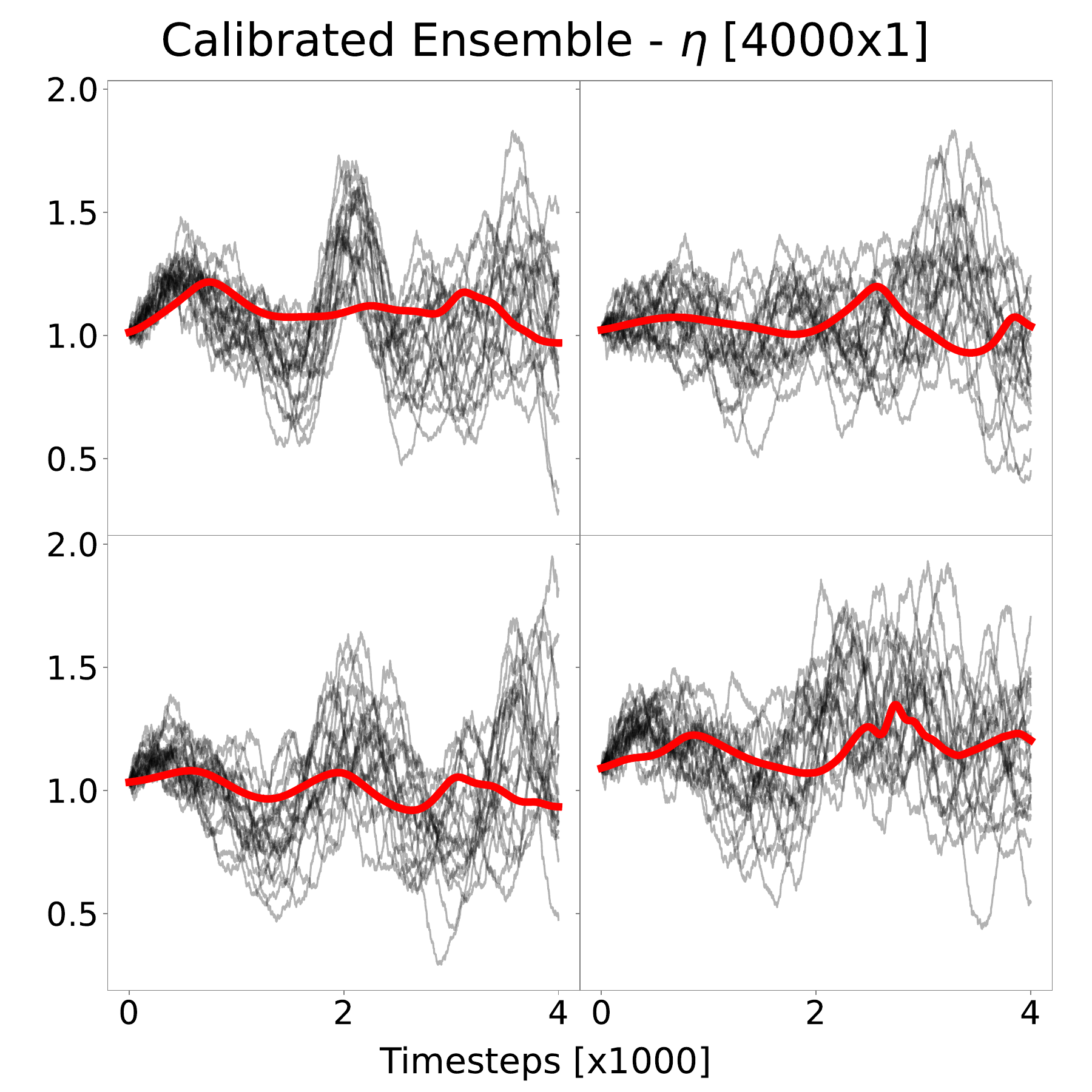}
        \caption{Height~$\eta$}
    \end{subfigure}
    \begin{subfigure}{0.3\linewidth}
        \includegraphics[width=\linewidth]{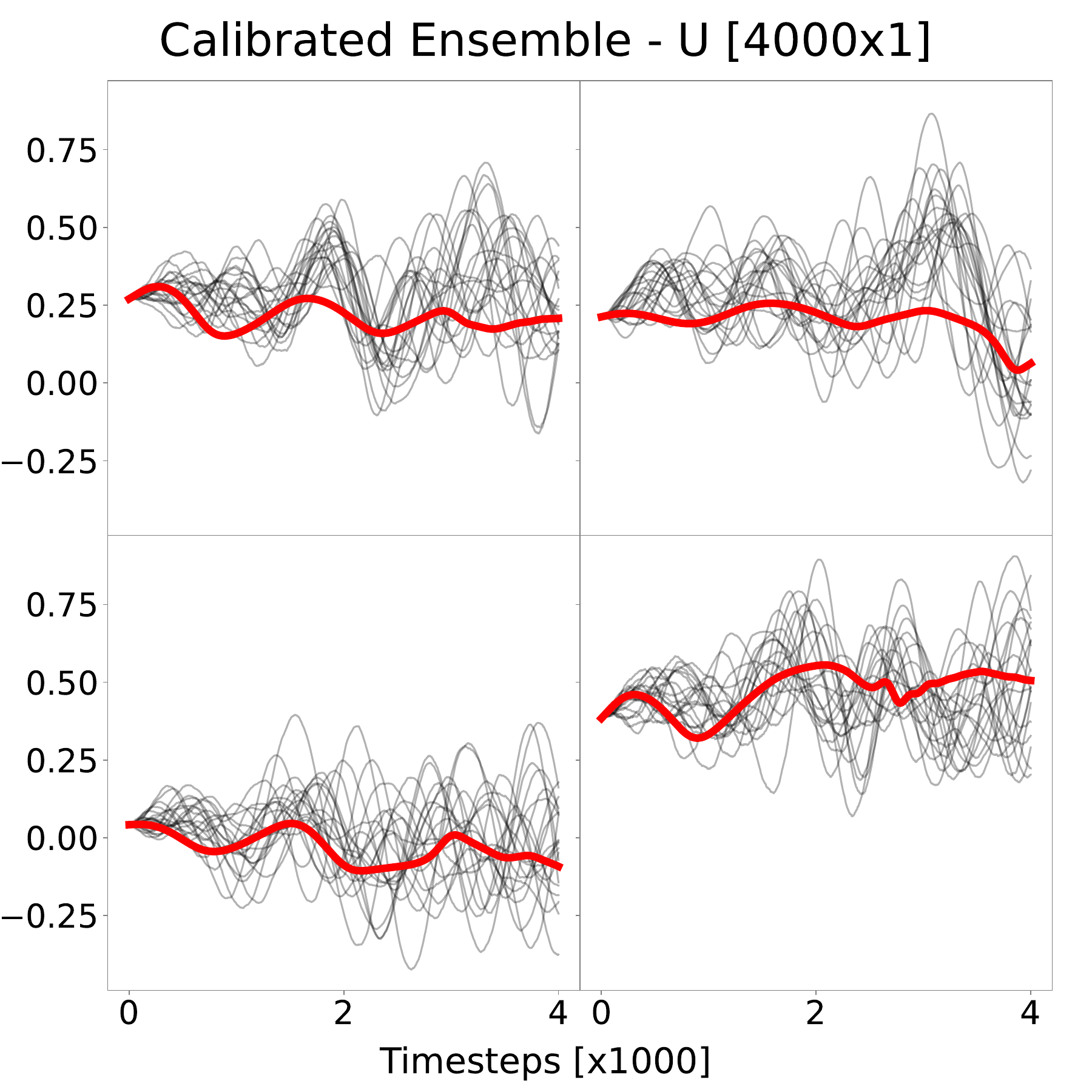}
        \caption{Velocity U}
    \end{subfigure}
    \begin{subfigure}{0.3\linewidth}
        \includegraphics[width=\linewidth]{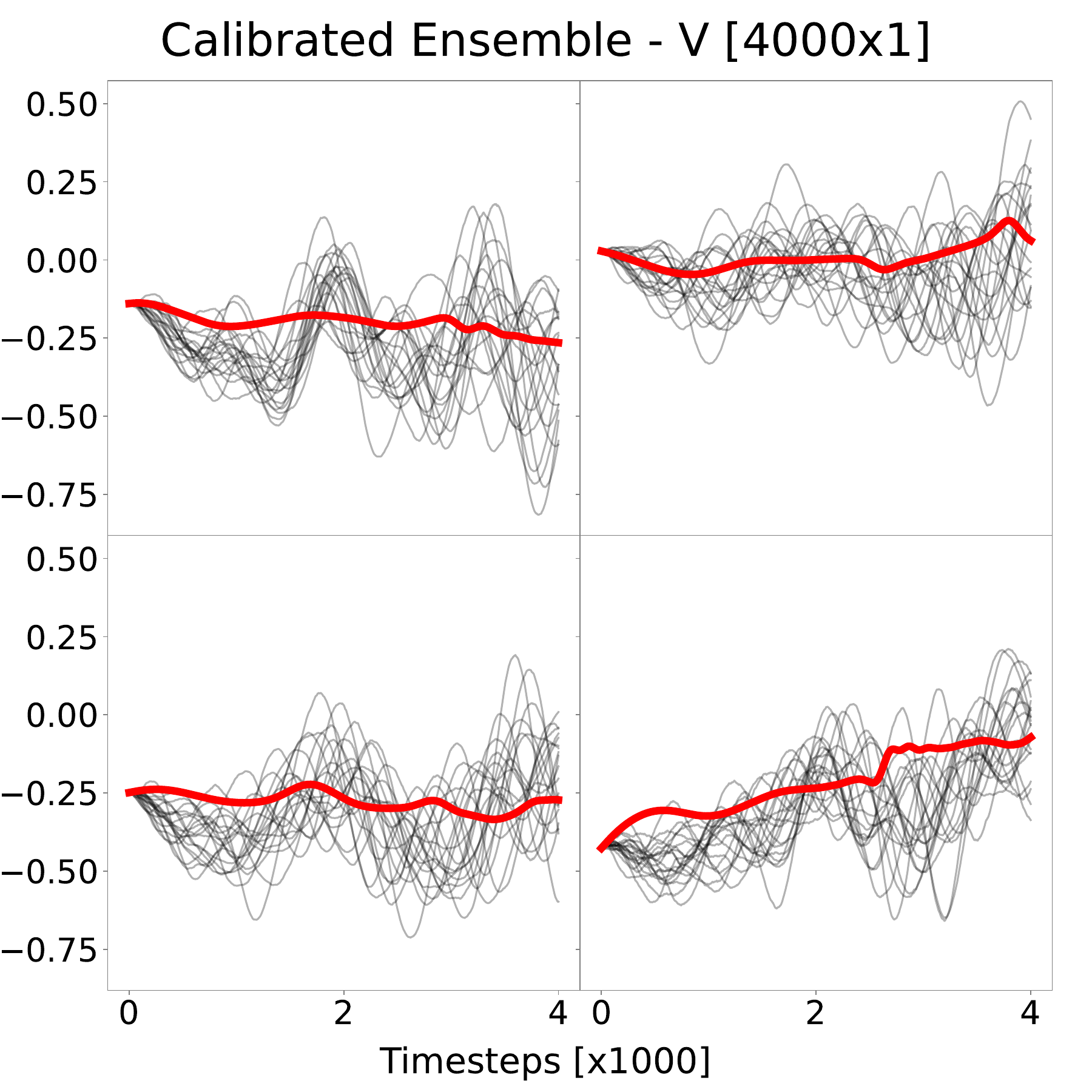}
        \caption{Velocity V}
    \end{subfigure}
    \caption{Ensemble evolutions (black lines) for forecast ensembles of the shallow water variables $\eta$ (left), $U$ (center), and $V$ (right) at four different sample grid points. The ensemble forecast evolutions are shown for a forecast horizon of $4000$ timesteps (horizontal axes) of size $\delta t =0.0001$ and no initial ensemble perturbation. The ensemble consists of $N_{\text{ens}}=20$ members. Red lines show the evolution of the PDE on the fine grid.}
    \label{fig:calib_spaghetti}
\end{figure}

Finally, we inspect the evolution of the forecast ensembles by means of their trajectories at specified sample grid points around the center of the domain. These are plotted for all three variables with a forecast horizon of $4000$ steps in Figure~\ref{fig:calib_spaghetti} together with the trajectories of the fine-grained truth. We observe that for all variables the ensembles spread around the truth.
\subsection{Data Assimilation}\label{sec:num_da}

Data assimilation is performed using a particle filter with tempering and jittering for the SRSW model variable $\eta$ with trained noise generator. The filter is assessed in several different configurations. We chose to run the ensemble for a total of 400 timesteps with different assimilation windows of $N_f = 10$ and $N_f=20$ timesteps. Different ensemble sizes of $N_{ens}=50$, $100$ and $200$ particles are used and we simultaneously observe a noisy version of the true fine-grid trajectory of the deterministic model at $d_{obs}=1$, $4$, $9$ and $16$ grid locations, respectively. The observation noise is set to Gaussian with a standard deviation of $\sigma_{\text{obs}}=0.01$. As error metrics we report, at assimilation times, the ensemble bias
\begin{equation}
    \text{Bias}^{i,j}(t) = \left(\frac{1}{N_{ens}} \sum_{p=1}^{N_{ens}} x^{i,j}_p(t)\right) - x^{i,j}_{\text{true}}(t)
\end{equation}
and the ensemble root mean square error
\begin{equation}
    \text{RMSE}^{i,j}(t) = \sqrt{\frac{1}{N_{ens}} \sum_{p=1}^{N_{ens}} [x^{i,j}_p(t) - x^{i,j}_{\text{true}}(t)]^2}
\end{equation}
with respect to the deterministic, true trajectory at an observed grid location $(i,j)$.

\begin{table}[ht!]\small
    \centering
    \caption{\textbf{Bias.} Time average of ensemble bias with respect to the truth. Reported is the best of all observed locations.}
    \begin{tabular}{c c c c c c c c c}
    \toprule
          \multirow[c]{2}*{$\mathbf{N_{\text{ens}}}$} & \multicolumn{2}{c}{$\mathbf{d_{\text{obs}}=1}$} & \multicolumn{2}{c}{$\mathbf{d_{\text{obs}}=4}$} & \multicolumn{2}{c}{$\mathbf{d_{\text{obs}}=9}$} & \multicolumn{2}{c}{$\mathbf{d_{\text{obs}}=16}$} \\
        \cmidrule(lr){2-3} \cmidrule(lr){4-5} \cmidrule(lr){6-7} \cmidrule(lr){8-9}
            & $\mathbf{N_{f}=10}$ & $\mathbf{N_{f}=20}$ & $\mathbf{N_{f}=10}$ & $\mathbf{N_{f}=20}$ & $\mathbf{N_{f}=10}$ & $\mathbf{N_{f}=20}$ & $\mathbf{N_{f}=10}$ & $\mathbf{N_{f}=20}$\\
        \midrule
         $\mathbf{50}$ & $0.0070$ & \cellcolor{black!15!white}$0.0048$ &$0.0070$& $0.0065$ &$0.0099$& $0.0181$&         $0.0097$    & $0.0118$\\
         $\mathbf{100}$ & $0.0059$ & $0.0059$ & $0.0061$ & $0.0069$ & \cellcolor{black!15!white}$0.0088$ & $0.0117$ &   $0.0097$    & $0.0121$\\
         $\mathbf{200}$ & $0.0061 $& $0.0048$ & $0.0066$ & \cellcolor{black!15!white}$0.0061$ & $0.0096$ & $0.0109$&    \cellcolor{black!15!white}$0.0081$    & $0.0144$\\
        \bottomrule
    \end{tabular}
    \label{tab:bias}
\end{table}

The smallest mean bias over time measured at the observed grid locations for each scenario is shown in Table~\ref{tab:bias}. The best values for a given observation dimension $d_{\text{obs}}$ are highlighted. We observe that the filter performs well in all cases analyzed. The bias remains well below the observation noise. For the case when only one observation is assimilated there are no significant improvements when the number of particles is increased, but the results are better when the assimilation windows are halved. For the case when four observations are assimilated there are improvements both as the number of particles is increased and the assimilation windows are halved. No distinguishable pattern can be observed when 9, respectively, 16 observations are assimilated.

\begin{table}[ht!]\small
    \centering
    \caption{\textbf{RMSE.} Time average of ensemble RMSE with respect to the truth. Reported is the best of all observed locations.}
    \begin{tabular}{c c c c c c c c c}
    \toprule
          \multirow[c]{2}*{$\mathbf{N_{\text{ens}}}$} & \multicolumn{2}{c}{$\mathbf{d_{\text{obs}}=1}$} & \multicolumn{2}{c}{$\mathbf{d_{\text{obs}}=4}$} & \multicolumn{2}{c}{$\mathbf{d_{\text{obs}}=9}$} & \multicolumn{2}{c}{$\mathbf{d_{\text{obs}}=16}$} \\
        \cmidrule(lr){2-3} \cmidrule(lr){4-5} \cmidrule(lr){6-7} \cmidrule(lr){8-9}
            & $\mathbf{N_{f}=10}$ & $\mathbf{N_{f}=20}$ & $\mathbf{N_{f}=10}$ & $\mathbf{N_{f}=20}$ & $\mathbf{N_{f}=10}$ & $\mathbf{N_{f}=20}$ & $\mathbf{N_{f}=10}$ & $\mathbf{N_{f}=20}$\\
        \midrule
         $\mathbf{50}$ & $0.0170$ & $0.0188$ &\cellcolor{black!15!white}$0.0160$& $0.0195$ &\cellcolor{black!15!white}$0.0161$& $0.0256$&$0.0161$& $0.0194$\\
         $\mathbf{100}$ & \cellcolor{black!15!white}$0.0168$ & $0.0202$ & $0.0164$ & $0.0191$&$0.0161$& $0.0214$                        &$0.0161$& $0.0206$\\
         $\mathbf{200}$ & $0.0169$& $0.0197$ & $0.0170$ & $0.0199$&$0.0172$& $0.0212$                                                   &\cellcolor{black!15!white}$0.0142$& $0.0216$\\
        \bottomrule
    \end{tabular}
    \label{tab:rmse}
\end{table}

The smallest mean RMSE over time among the observed grid locations for each scenario is shown in Table~\ref{tab:rmse}. The best values for a given observation dimension $d_{\text{obs}}$ are highlighted. We observe that the filter performs reasonably well also insofar as the RMSE is concerned. However no improvements are observed when the number of particles is increased or the assimilation intervals are halved.

\begin{figure}
    \centering
    \begin{subfigure}{0.8\textwidth}
        \includegraphics[width=\linewidth]{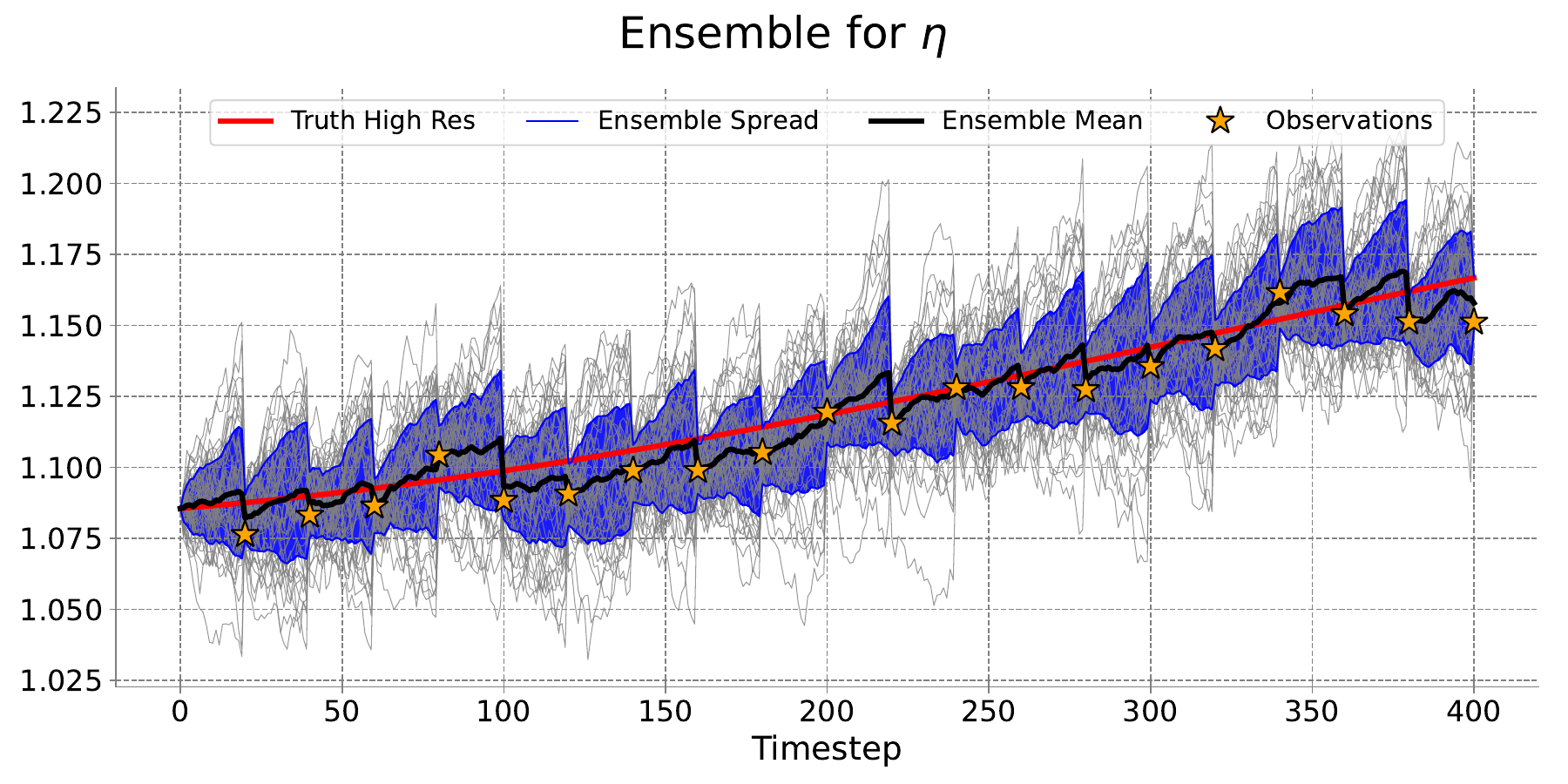}
        \caption{Ensemble Evolution}
        \label{fig:filt1d-ensemble}
    \end{subfigure}
    \hfill
    \begin{subfigure}{0.4\textwidth}
        \includegraphics[width=\linewidth]{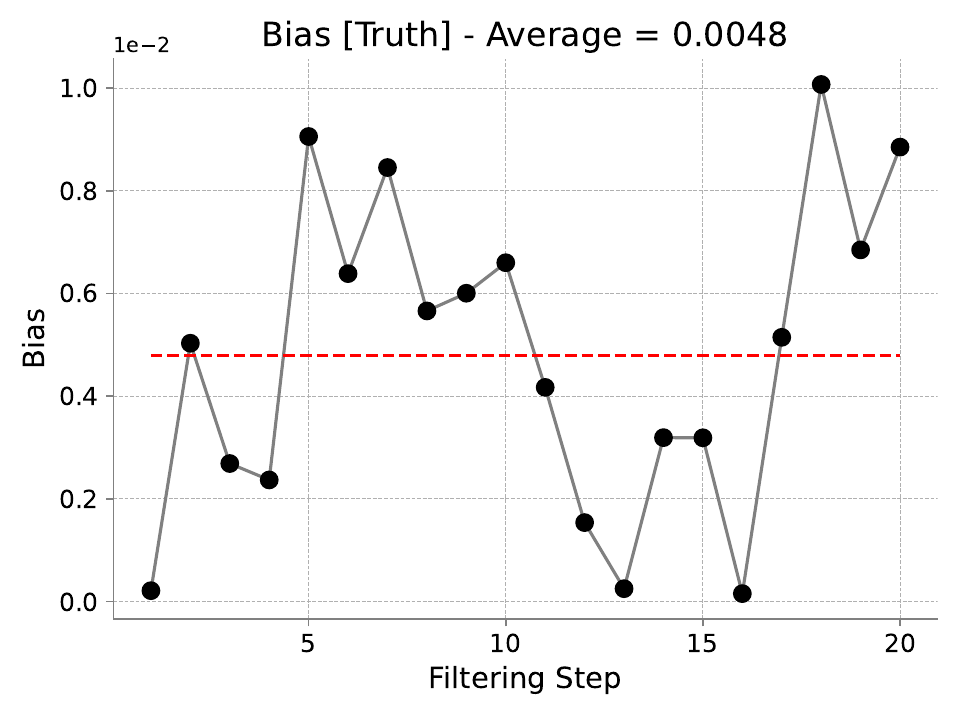}
        \caption{Ensemble Bias}
        \label{fig:filt1d-bias}
    \end{subfigure}
    \begin{subfigure}{0.4\textwidth}
        \includegraphics[width=\linewidth]{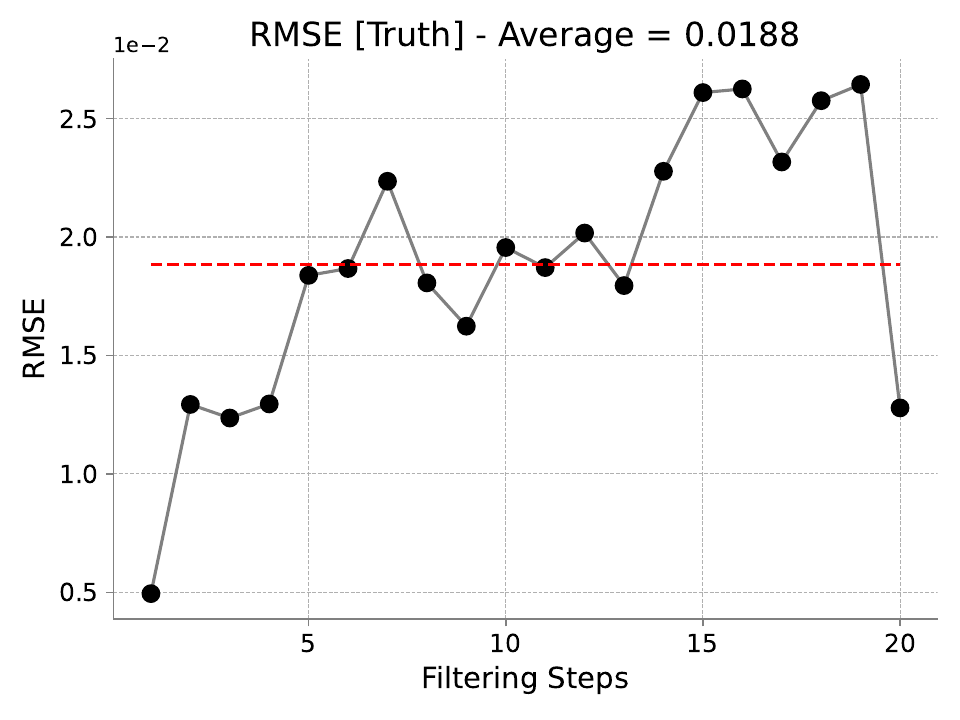}
        \caption{Ensemble RMSE}
        \label{fig:filt1d-rmse}
    \end{subfigure}
    \caption{\textbf{Filtering Results 1D.} Results of the filtering experiment over $400$ total timesteps. The assimilated variable is $\eta$ using an ensemble of $50$ particles and an assimilation window of $20$ forecast timesteps. We show, at the observed gird location, (a) the ensemble evolution compared to the true deterministic fine-grid trajectory, (b) the ensemble bias over time, and (c) the ensemble RMSE.}
    \label{sp:1Dobs}
\end{figure}

\begin{figure}
    \centering
    \begin{subfigure}{0.8\textwidth}
        \includegraphics[width=\linewidth]{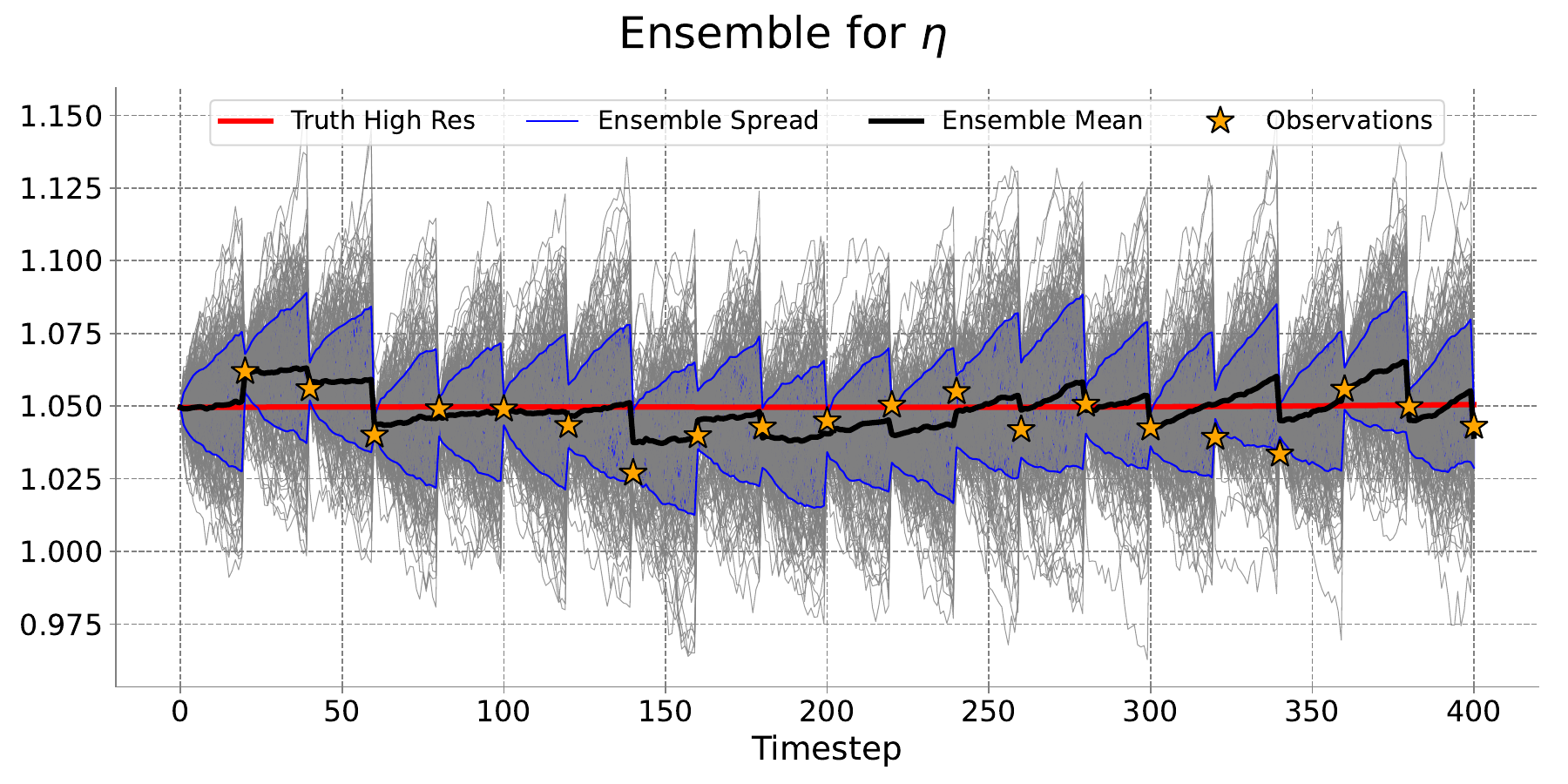}
        \caption{Ensemble Evolution}
        \label{fig:filt4d-ensemble}
    \end{subfigure}
    \hfill
    \begin{subfigure}{0.4\textwidth}
        \includegraphics[width=\linewidth]{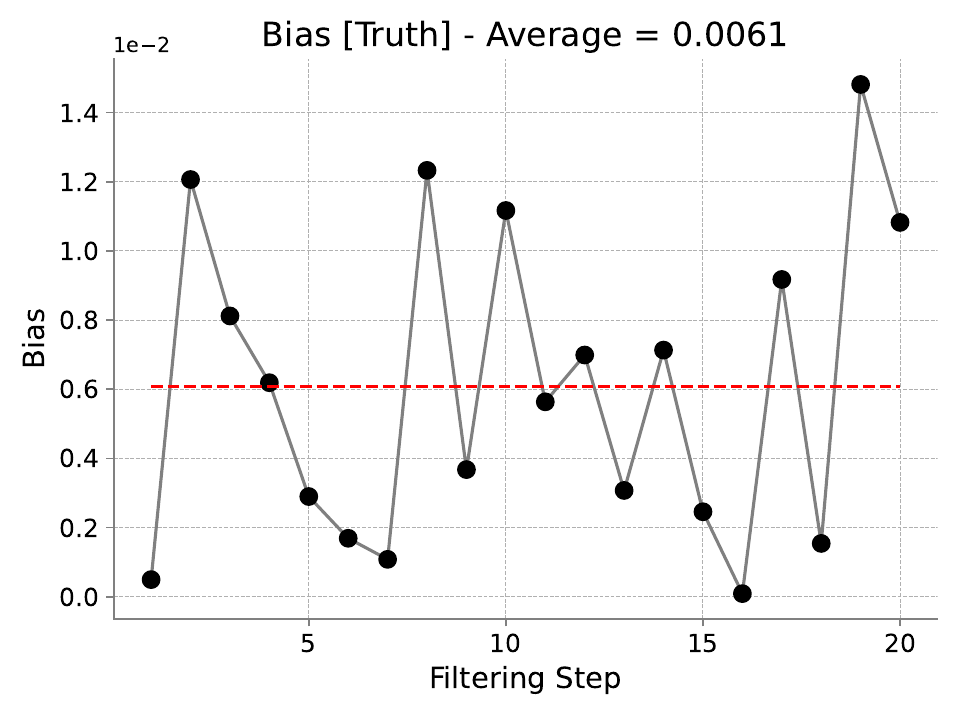}
        \caption{Ensemble Bias}
        \label{fig:filt4d-bias}
    \end{subfigure}
    \begin{subfigure}{0.4\textwidth}
        \includegraphics[width=\linewidth]{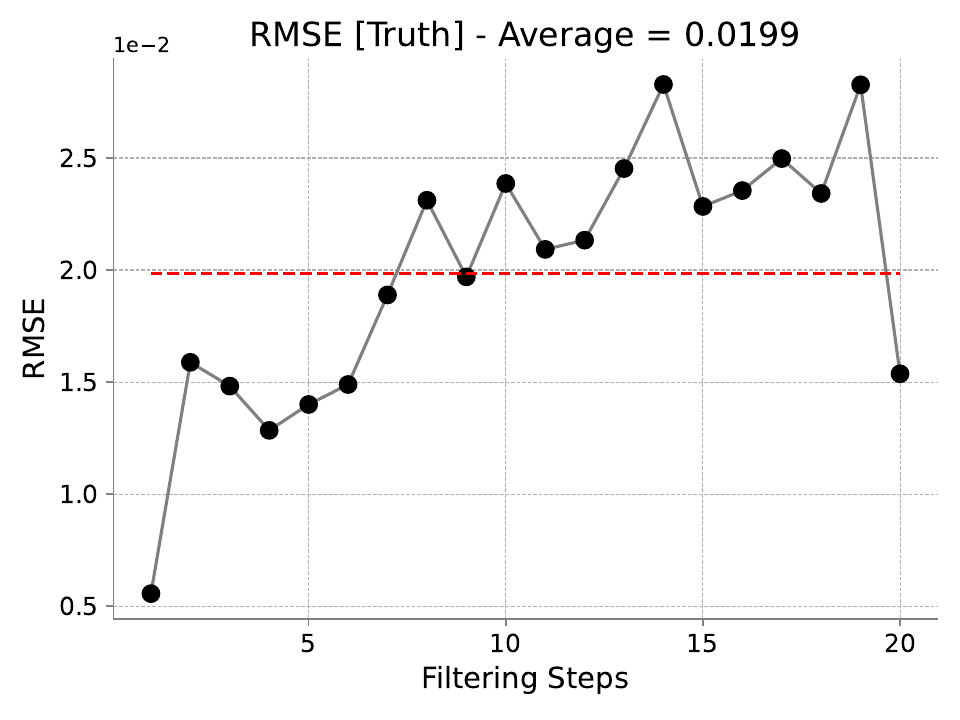}
        \caption{Ensemble RMSE}
        \label{fig:filt4d-rmse}
    \end{subfigure}
    \caption{\textbf{Filtering Results 4D.} Results of the filtering experiment over $400$ total timesteps. The assimilated variable is $\eta$ using an ensemble of $200$ particles and an assimilation window of $20$ forecast timesteps. We show, at the best observed gird location, (a) the ensemble evolution compared to the true deterministic fine-grid trajectory, (b) the ensemble bias over time, and (c) the ensemble RMSE.}
    \label{sp:4Dobs}
\end{figure}

\begin{figure}
    \centering
    \begin{subfigure}{0.8\textwidth}
        \includegraphics[width=\linewidth]{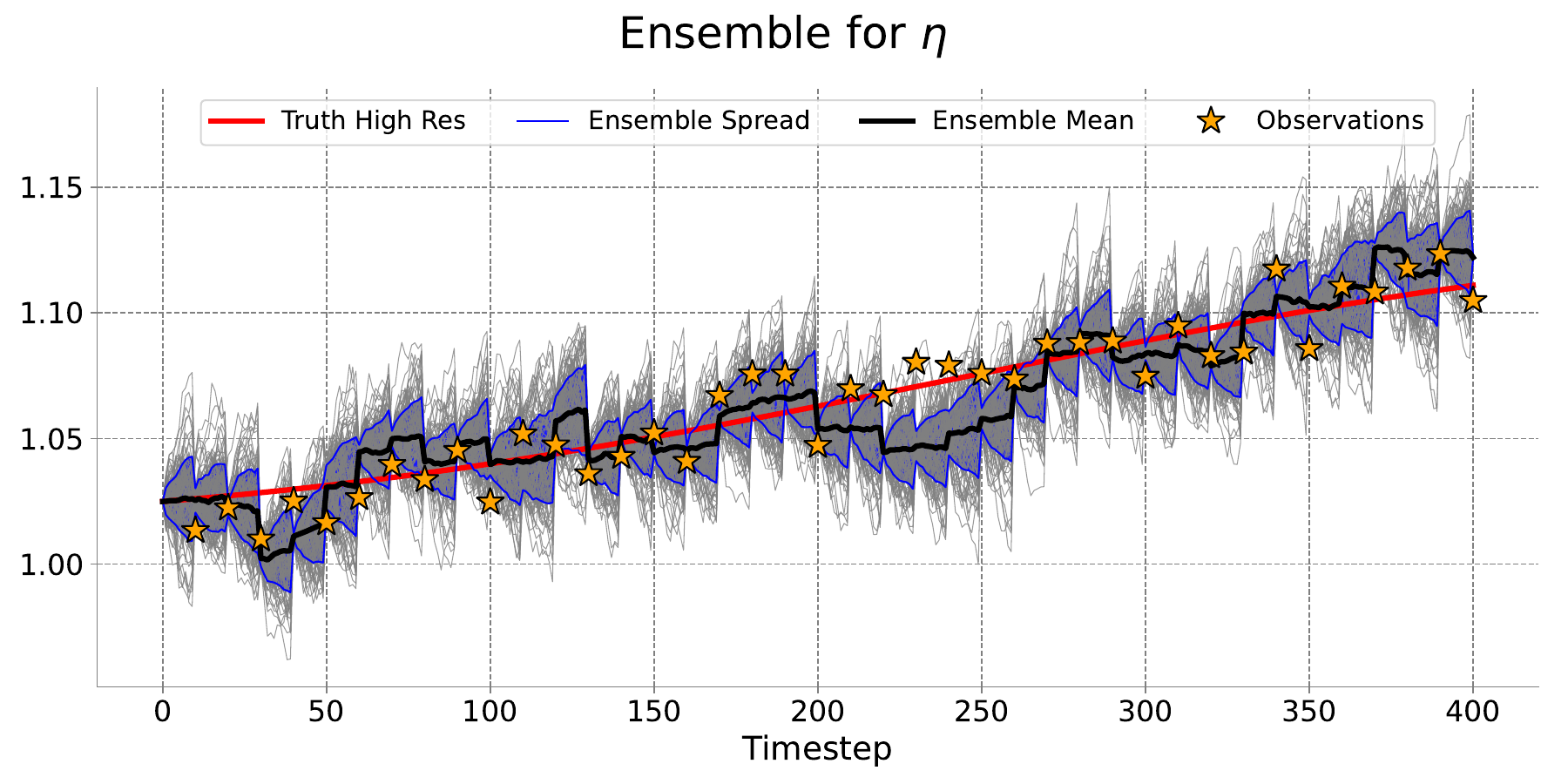}
        \caption{Ensemble Evolution}
        \label{fig:filt9d-ensemble}
    \end{subfigure}
    \hfill
    \begin{subfigure}{0.4\textwidth}
        \includegraphics[width=\linewidth]{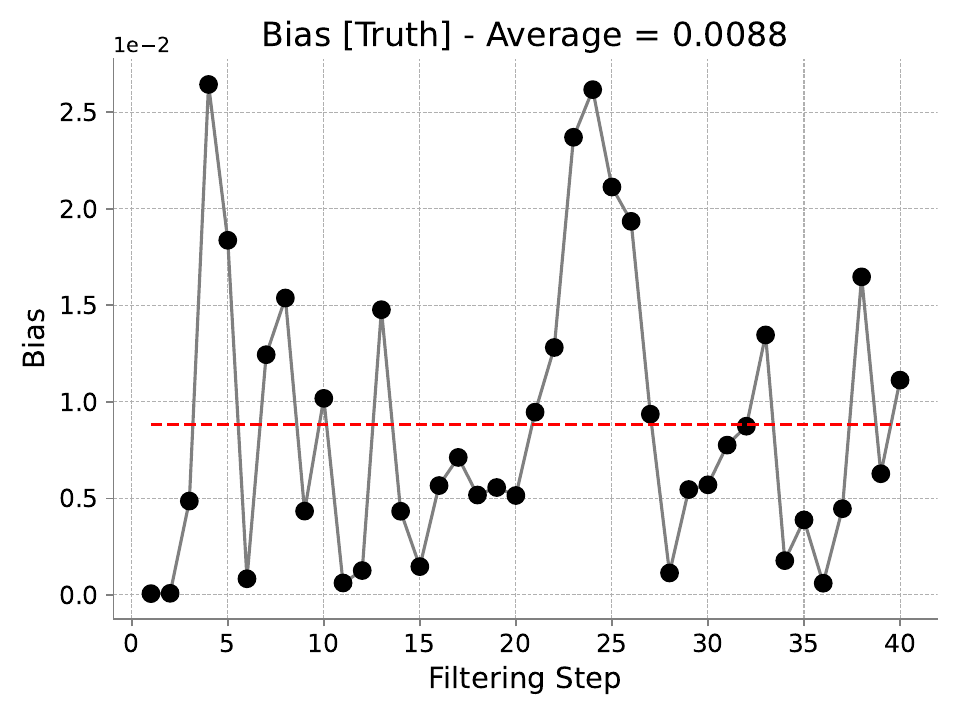}
        \caption{Ensemble Bias}
        \label{fig:filt9d-bias}
    \end{subfigure}
    \begin{subfigure}{0.4\textwidth}
        \includegraphics[width=\linewidth]{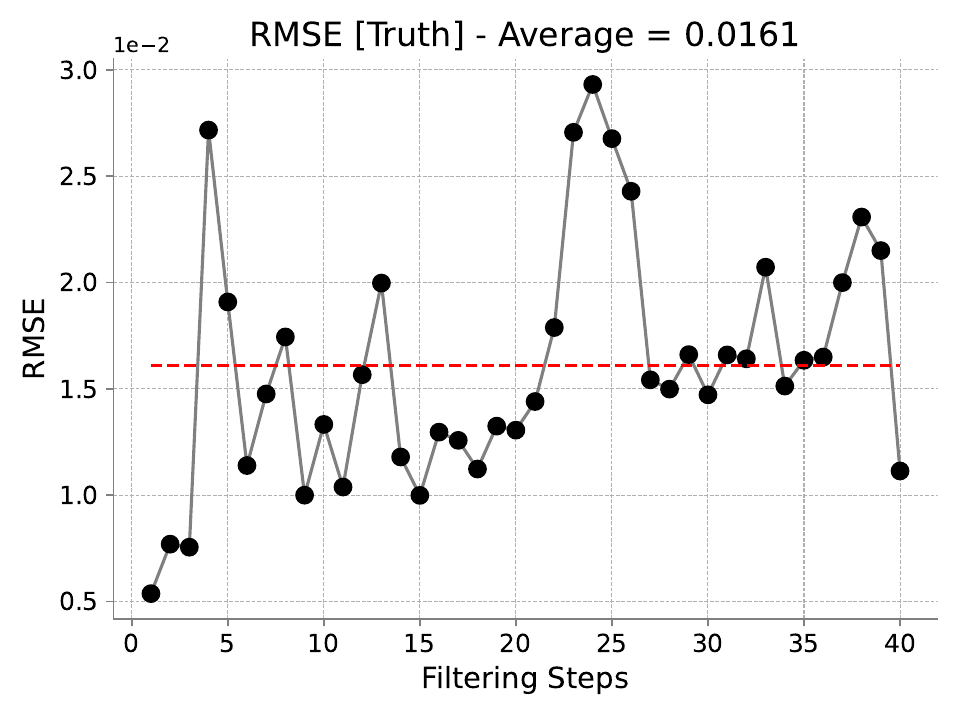}
        \caption{Ensemble RMSE}
        \label{fig:filt9d-rmse}
    \end{subfigure}
    \caption{\textbf{Filtering Results 9D.} Results of the filtering experiment over $400$ total timesteps. The assimilated variable is $\eta$ using an ensemble of $100$ particles and an assimilation window of $10$ forecast timesteps. We show, at the best observed gird location, (a) the ensemble evolution compared to the true deterministic fine-grid trajectory, (b) the ensemble bias over time, and (c) the ensemble RMSE.}
    \label{sp:9Dobs}
\end{figure}

\begin{figure}
    \centering
    \begin{subfigure}{0.8\textwidth}
        \includegraphics[width=\linewidth]{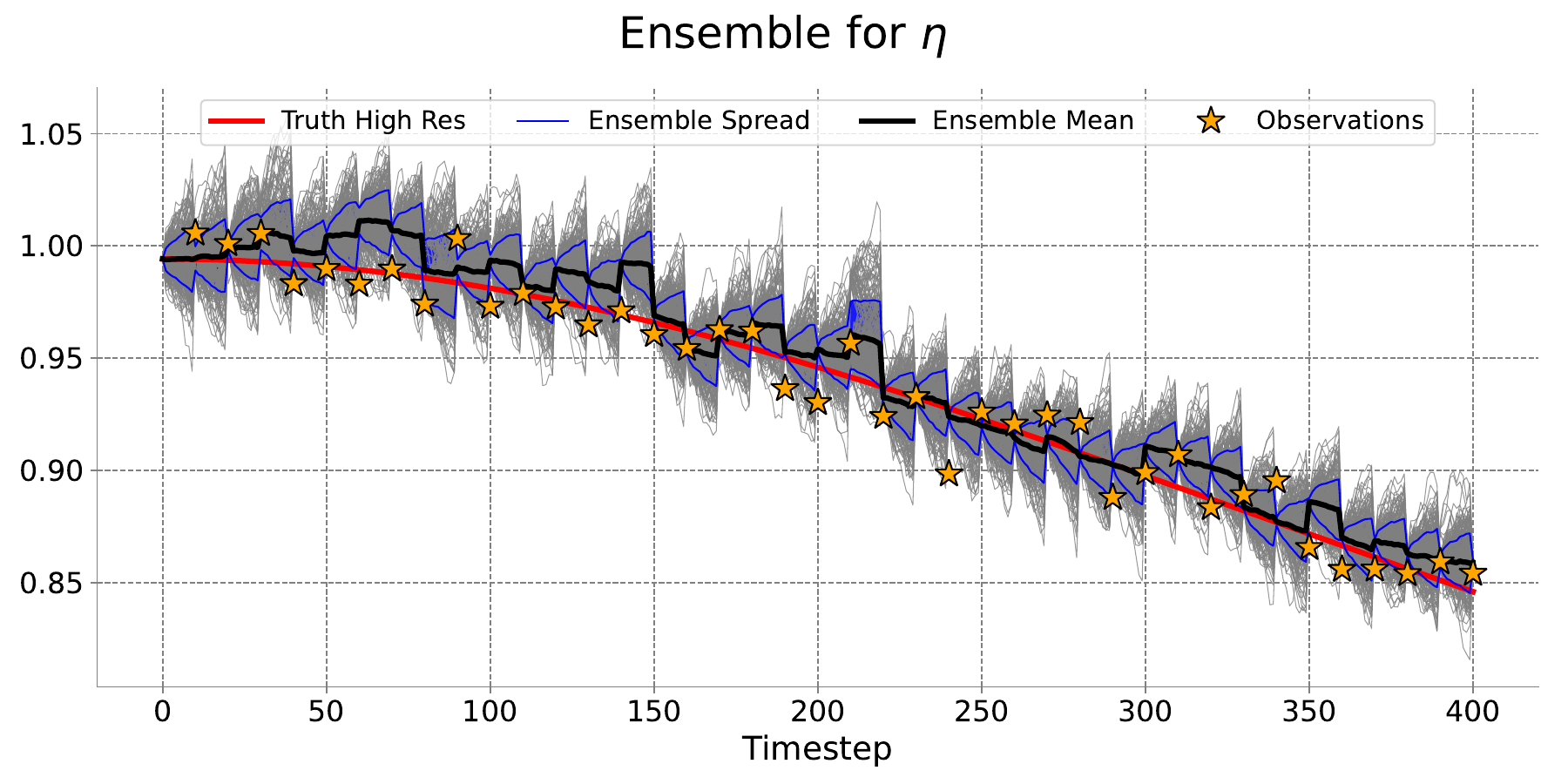}
        \caption{Ensemble Evolution}
        \label{fig:filt16d-ensemble}
    \end{subfigure}
    \hfill
    \begin{subfigure}{0.4\textwidth}
        \includegraphics[width=\linewidth]{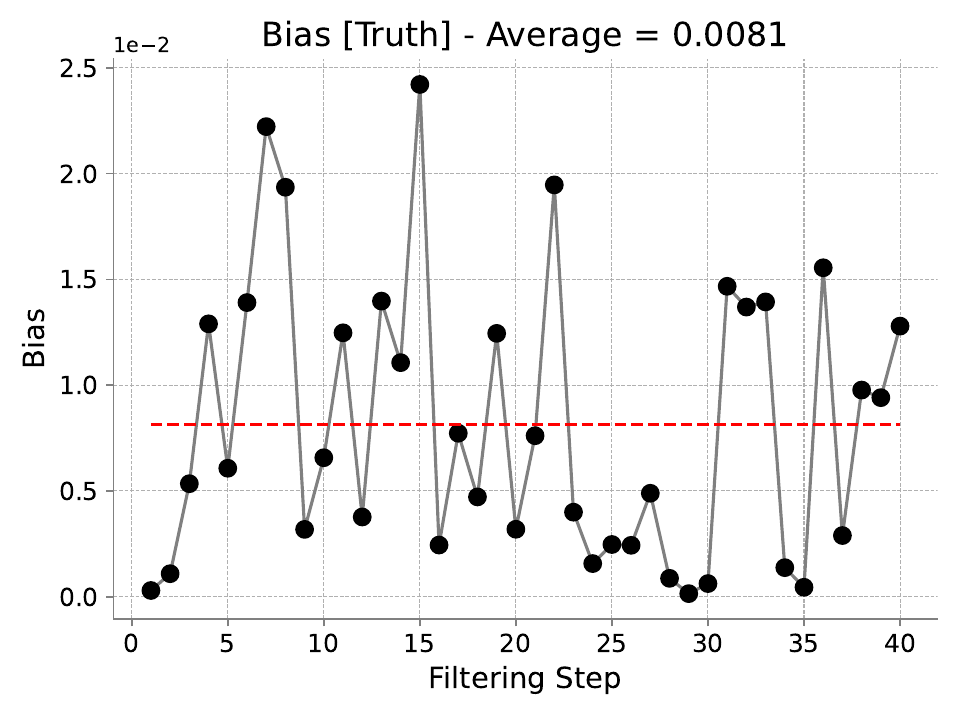}
        \caption{Ensemble Bias}
        \label{fig:filt16d-bias}
    \end{subfigure}
    \begin{subfigure}{0.4\textwidth}
        \includegraphics[width=\linewidth]{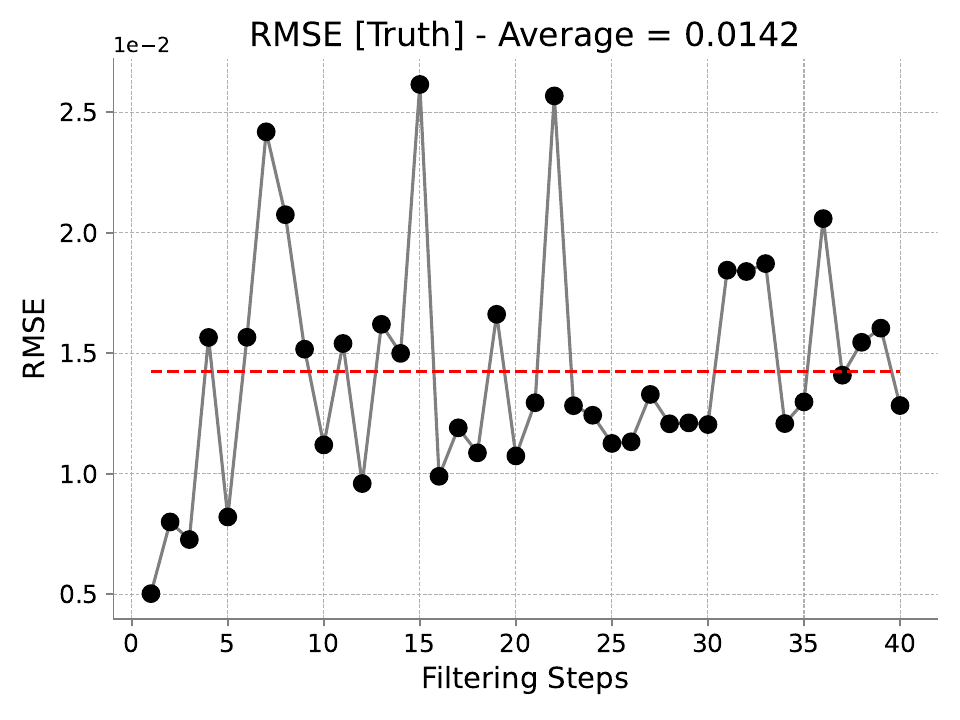}
        \caption{Ensemble RMSE}
        \label{fig:filt16d-rmse}
    \end{subfigure}
    \caption{\textbf{Filtering Results 16D.} Results of the filtering experiment over $400$ total timesteps. The assimilated variable is $\eta$ using an ensemble of $200$ particles and an assimilation window of $10$ forecast timesteps. We show, at the best observed gird location, (a) the ensemble evolution compared to the true deterministic fine-grid trajectory, (b) the ensemble bias over time, and (c) the ensemble RMSE.}
    \label{sp:16Dobs}
\end{figure}

\begin{figure}
    \centering
    \begin{subfigure}{0.8\textwidth}
        \includegraphics[width=\linewidth]{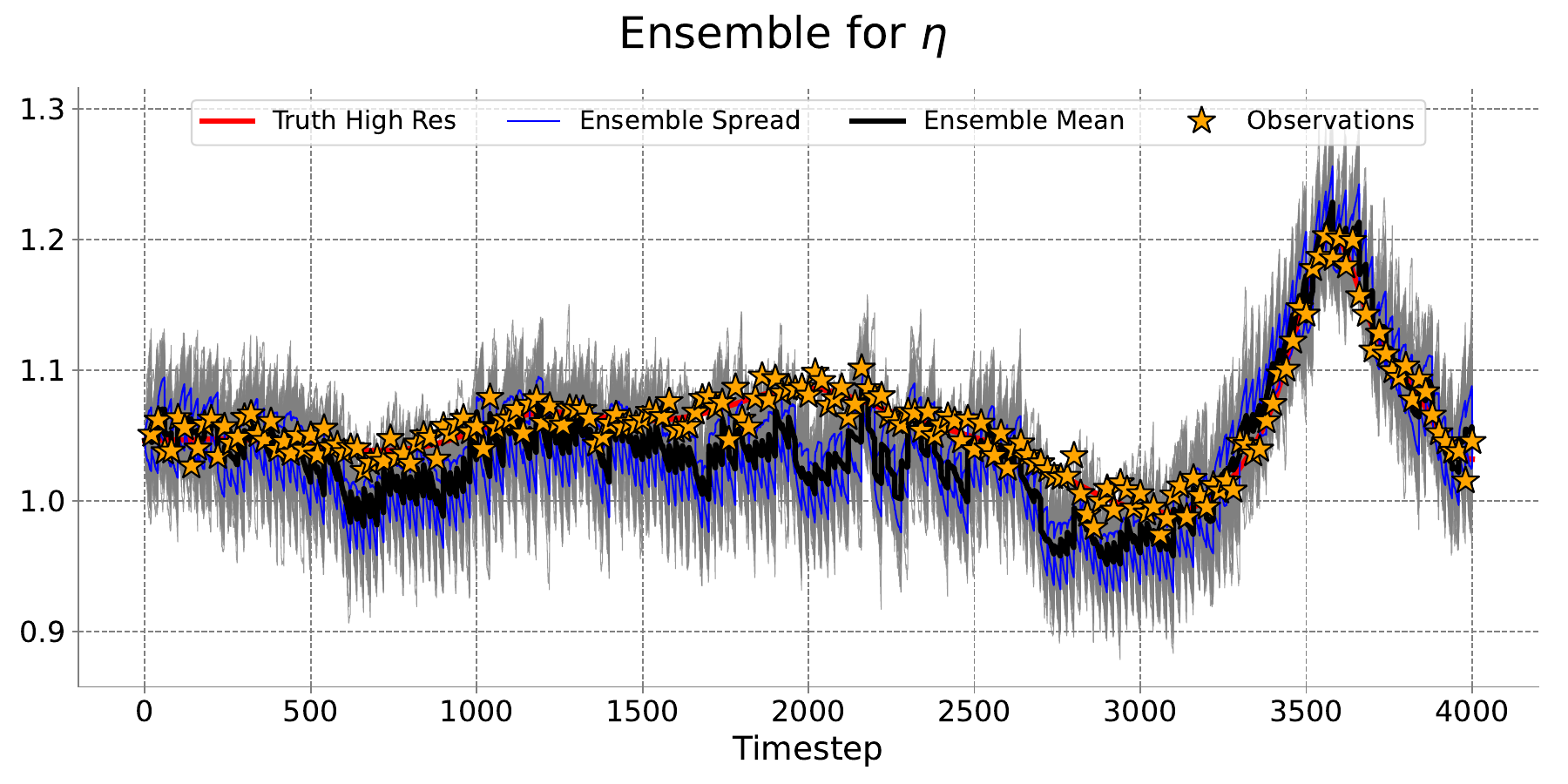}
        \caption{Ensemble Evolution}
        \label{fig:filt4d-long-ensemble}
    \end{subfigure}
    \hfill
    \begin{subfigure}{0.4\textwidth}
        \includegraphics[width=\linewidth]{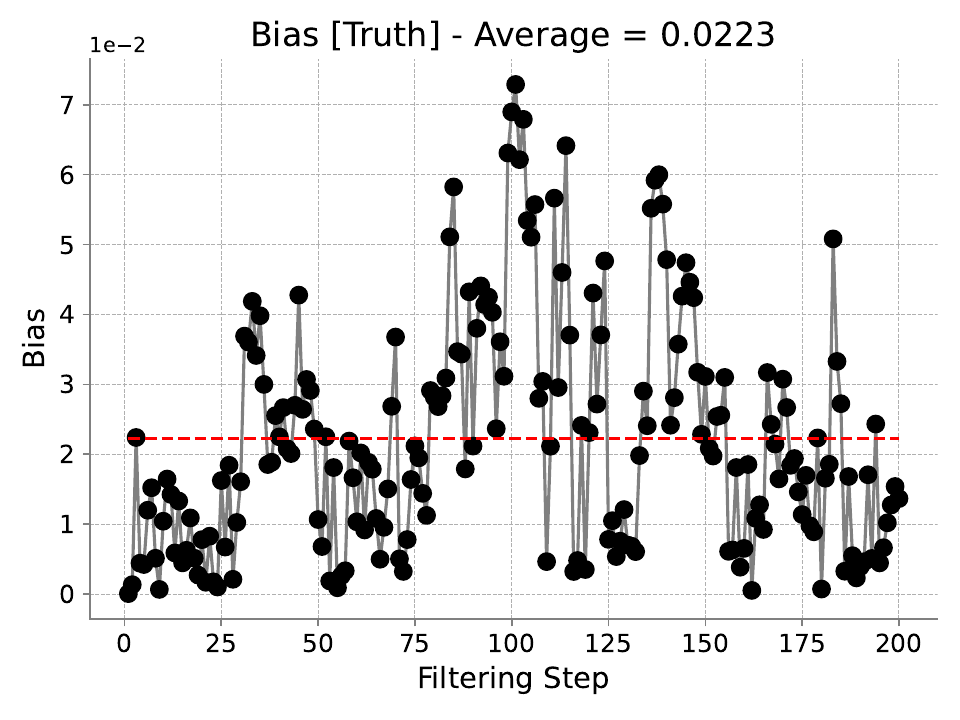}
        \caption{Ensemble Bias}
        \label{fig:filt4d-long-bias}
    \end{subfigure}
    \begin{subfigure}{0.4\textwidth}
        \includegraphics[width=\linewidth]{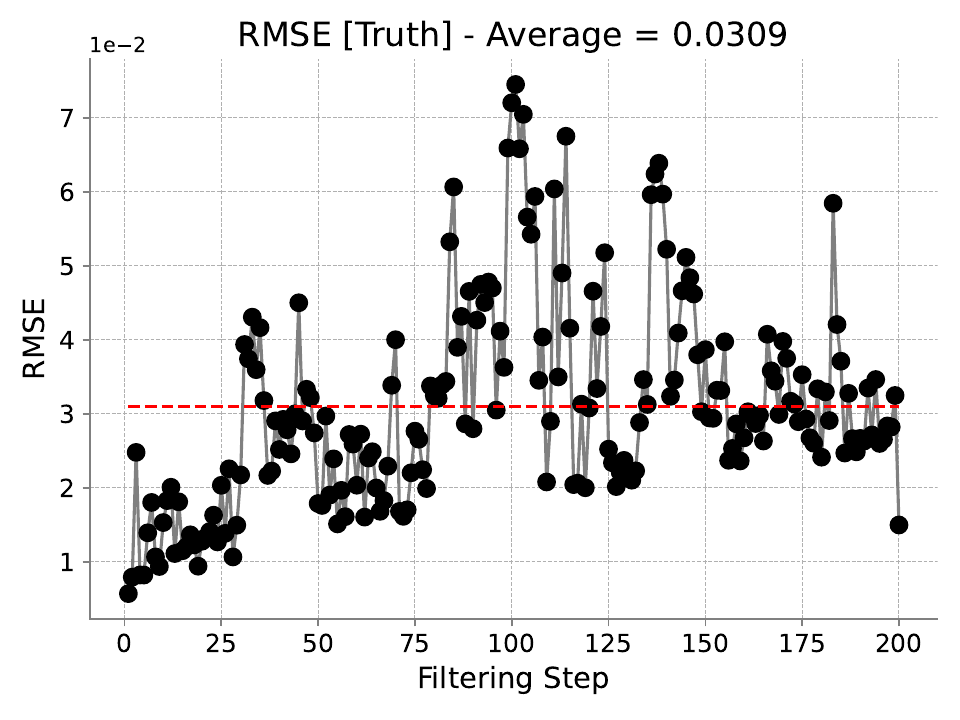}
        \caption{Ensemble RMSE}
        \label{fig:filt4d-long-rmse}
    \end{subfigure}
    \caption{\textbf{Long Filtering Results 4D.} Results of the filtering experiment over $4000$ total timesteps. The assimilated variable is $\eta$ using an ensemble of $200$ particles and an assimilation window of $20$ forecast timesteps. We show, at the best observed gird location, (a) the ensemble evolution compared to the true deterministic fine-grid trajectory, (b) the ensemble bias over time, and (c) the ensemble RMSE.}
    \label{sp:4Dobs:long}
\end{figure}

Figures \ref{sp:1Dobs}, \ref{sp:4Dobs}, and \ref{sp:9Dobs}, \ref{sp:16Dobs} represent the evolution of the height of the cloud of particles at a given point in the grid corresponding to the 1, 4, 9, and 16 observations assimilated, respectively. As expected the particles spread in between consecutive assimilation times and they are brought closer to the true trajectory each time new data is assimilated. A similar pattern is observed when the for the long run (4000 time-steps) presented in Figure \ref{sp:4Dobs:long}.

\section{Conclusions and future work}

In this paper we explore a novel calibration methodology using diffusion generative models combined with a particle filtering methodology augmented with tempering and jittering for dynamic state estimation. We generate synthetic data that statistically aligns with given set of observation (in this case the time increments of numerical approximations of a solution of the rotating shallow water equation). This allows us to efficiently implement a model reduction and assimilate data from a reference system state (the “truth'') modelled by a highly resolved numerical solution of the rotating shallow water equation on a spatial grid of size 128x128 into a stochastic system run on a spatial grid of size 32x32. The new samples are incorporated into a particle filtering methodology augmented with tempering and jittering for dynamic state estimation. This work demonstrates how generative models can be used to improve the predictive accuracy for particle filters, providing a more computationally efficient solution for data assimilation and model calibration. We enumerate below a number of improvements and directions to be pursued in future work.  

As explained in \cite{us} we use a specific type of noise to model the unresolved scales. In particular we assume that the noise is divergence free. This is very convenient as it enables us to calibrate the noise only by using data coming from recording the $\eta$ the height of the fluid column. This ansatz may not be justified in general. The constraint can be lifted at the expense of using additional observational data, for example by observing one co-ordinate of the velocity field. This will possibly improve the accuracy of the particle filter. In this study we calibrate the noise using data coming from $\eta$ the height variable and performe the particle filter on the same variable, However, it would be useful to incorporate the velocity field into the analysis and enlarge the study to incorporate bayesian inference of the associated velocity field.

As we explained, the particle filter gives reasonable results despite the fact that it is run on a coarser scale than that of the true signal. This is made possible by exploiting the robustness of the filter with respect to the transition model misspecification and using a judicious calibration of the unresolved scales.  However, we can see that the performance cannot be visibly improved by increasing the number of particles or reducing the size of the data assimilation intervals. The resulting systematic error is higher than the statistical error. To improve this result, we intend to combine the particle filter with an additional procedure called \emph{nudging}. Nudging is a popular algorithmic strategy in numerical filtering to deal with the problem of inference in high-dimensional dynamical systems. We will use it to tackle the transition model misspecification. Specifically, we will rely on the formulation of the nudging as a general operation increasing the likelihood function and reduce the systematic errors.

In future work, we will test this newly developed methodology to real data.  In particular we will assimilate sea surface height from a region in the Pacific that straddles the Equator from 37.7°S to 24.9°N and 176°E to 117.5°W. The data used for calibrating the stochastic model will be monthly sea-surface height at 0.1° resolution, from 1984 to 2016 (372 months). In this case the stochastic model will be the rotating shallow water equation with a boundary condition given by the projection of the sea surface height data onto a coarser resolution grid.  The data comes from Ocean Reanalysis System 5 (ORAS5)~\cite{oras5-data}. 

\begin{figure}
    \centering
    \begin{subfigure}{0.38\linewidth}
        \includegraphics[width=\linewidth]{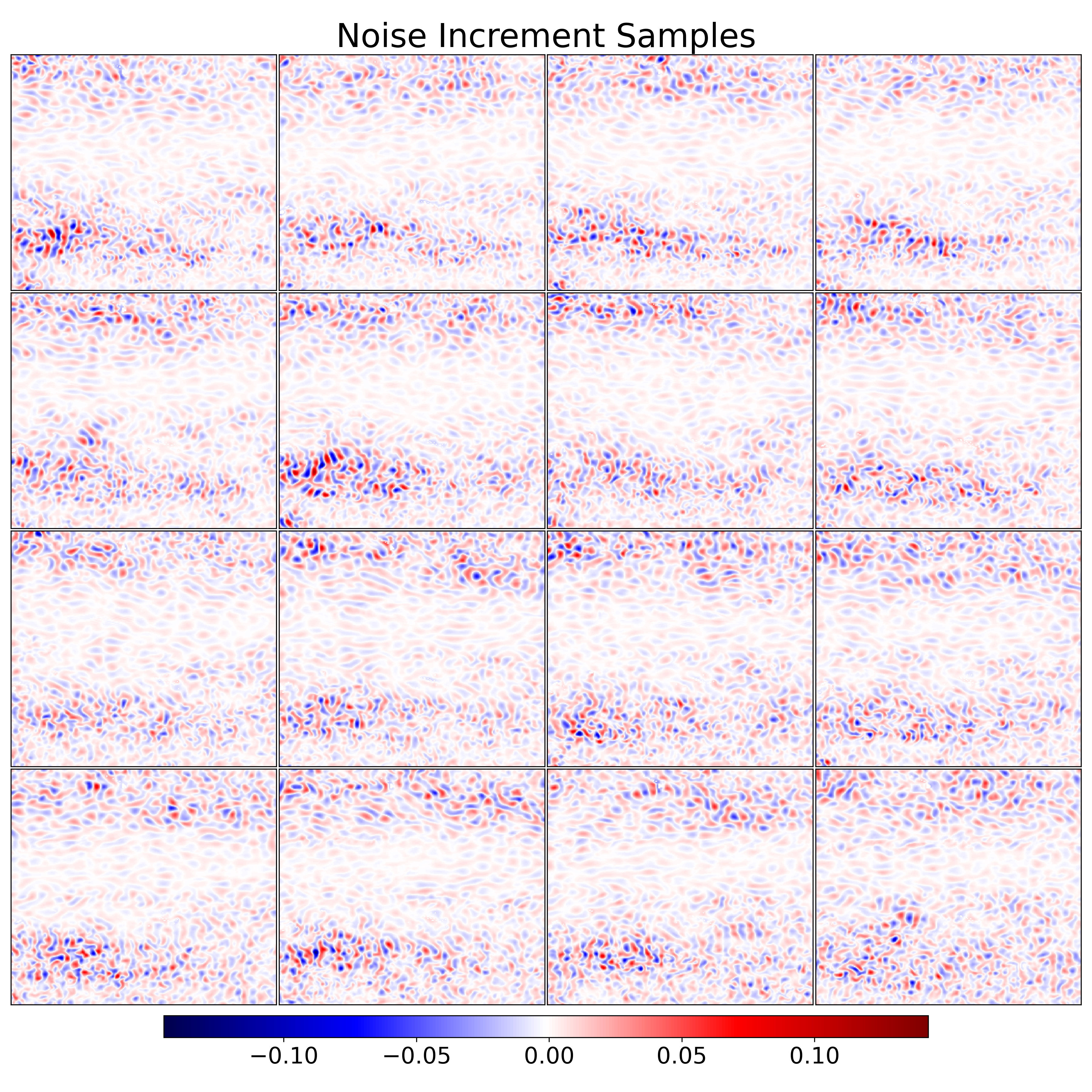}
        \caption{Noise Increments.}
        \label{fig:real_noise_samples}
    \end{subfigure}
    \begin{subfigure}{0.52\linewidth}
        \includegraphics[width=\linewidth]{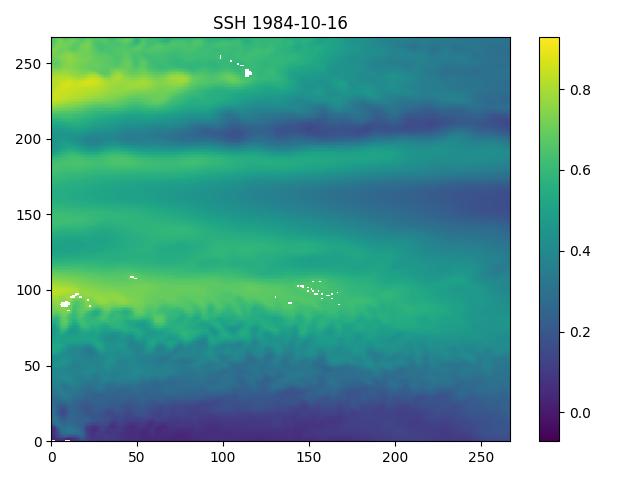}
        \caption{Sample of the Sea Surface height.}
        \label{fig:real_ssh}
    \end{subfigure}
    \caption{We obtained real monthly measurements of sea-surface height (SSH) from ORAS5. The measurements are taken from a square region in the pacific ocean centered at Lon: -150.732422, Lat: -7.449624. (a) Noise increments have been computed from the data and a sample of these is shown. The noise increments shown correspond to the input into the calibration equation. (b) A measurement of SSH from October 1984, taken from our ORAS5 dataset.}
    \label{fig:real-data}
\end{figure}

In Figure \ref{fig:real-data} we present 16 sample of the noise increment and one instance of the sea surface height from October 1984.  The area is in the Pacific and straddles the Equator from 37.7°S to 24.9°N and 176°E to 117.5°W. The noise increments show the signature of mesoscale eddies away from the equatorial band. \footnote{We thank Eviatar Bach for providing the data ocean sea surface height and the corresponding interpretation the  noise increments.}

\vspace{5mm}
\noindent\textbf{Funding}\\
\noindent All three authors have been supported by the European Research Council (ERC) under the European Union’s Horizon 2020 Research and Innovation Programme (ERC, Grant Agreement No 856408). 

\vspace{3mm}
\noindent\textbf{Data availability statement} \\
\noindent Data sharing not applicable to this article as no datasets were generated or analysed during the current study.

\vspace{3mm}
\noindent\textbf{Conflict of interest statement} \\
On behalf of all authors, the corresponding author states that there is no conflict of interest.

\clearpage
\bibliographystyle{plain}
\bibliography{main}

\end{document}